\documentclass[12pt,reqno]{amsart}

\makeatletter
\def\@settitle{\begin{center}%
  \baselineskip14\p@\relax
  \bfseries
  \@title
  \end{center}%
}

\makeatother

\usepackage{amsmath,amsfonts,amsthm,amsopn,amssymb,extarrows,multirow,footnote,float}
\usepackage{cite,marginnote}
\usepackage{booktabs}
\usepackage{bm}
\usepackage{color,enumitem,graphicx}
\usepackage[colorlinks=true,urlcolor=blue,
citecolor=red,linkcolor=blue,linktocpage,pdfpagelabels,
bookmarksnumbered,bookmarksopen]{hyperref}
\usepackage[english]{babel}

\usepackage[left=2.9cm,right=2.9cm,top=2.8cm,bottom=2.8cm]{geometry}
\usepackage[hyperpageref]{backref}

\numberwithin{equation}{section}

\makeindex
\makeindex
\def\lb{\lambda}

\def\N{\mathbb{N}}

\def\mS{\mathcal{S}}

\newtheorem{theorem}{Theorem}[section]
\newtheorem{definition}[theorem]{Definition}
\newtheorem{lemma}[theorem]{Lemma}
\newtheorem{corollary}[theorem]{Corollary}
\newtheorem{proposition}[theorem]{Proposition}
\newtheorem{remark}[theorem]{Remark}

\newcommand{\s}{\section}

\newcommand{\R}{\mathbb R}

\newcommand{\lab}{\label}
\newcommand{\bt}{\begin{theorem}}
	\newcommand{\et}{\end{theorem}}
\newcommand{\bl}{\begin{lemma}}
	\newcommand{\el}{\end{lemma}}
\newcommand{\bd}{\begin{definition}}
	\newcommand{\ed}{\end{definition}}
\newcommand{\bc}{\begin{corollary}}
	\newcommand{\ec}{\end{corollary}}
\newcommand{\bp}{\begin{proof}}
	\newcommand{\ep}{\end{proof}}
\newcommand{\bx}{\begin{example}}
	\newcommand{\ex}{\end{example}}
\newcommand{\bi}{\begin{exercise}}
	\newcommand{\ei}{\end{exercise}}
\newcommand{\bo}{\begin{proposition}}
	\newcommand{\eo}{\end{proposition}}
\newcommand{\br}{\begin{remark}}
	\newcommand{\er}{\end{remark}}
\newcommand{\beq}{\begin{equation}}
	\newcommand{\eeq}{\end{equation}}
\newcommand{\ba}{\begin{align}}
	\newcommand{\ea}{\end{align}}
\newcommand{\bn}{\begin{enumerate}}
	\newcommand{\en}{\end{enumerate}}
\newcommand{\bg}{\begin{align*}}
	\newcommand{\bcs}{\begin{cases}}
		\newcommand{\ecs}{\end{cases}}

	\newcommand{\bean}{\begin{eqnarray*}}
		\newcommand{\eean}{\end{eqnarray*}}

	\def\N{\mathbb{N}}

	\def\R{\mathbb{R}}

	\def\bd{\mathrm{bd}\,}

	\definecolor{my_red}{RGB}{255, 0, 0}

\def\@setauthors{\begin{center}%
  \def\and{\unskip\\}%
  \baselineskip16\p@\relax
  \normalfont % ¹Ø¼üÐÞ¸´£ºÇ¿ÖÆÊ¹ÓÃÕý³£×ÖÌå£¨±£Áô´óÐ¡Ð´£©
  \@author
  \end{center}%
}
	\renewcommand{\abstractname}{\normalfont\bfseries Abstract} % ¼Ó´ÖÇÒ·Ç´óÐ´
	\title[Normalized ground state for quasi-linear Schr\"odinger equation]{Serrin-type Overdetermination in Scaling Limits:  Sharp Existence and Asymptotic Behavior of Quasi-linear Schr\"{o}dinger Energy Ground States}
	
	\author[L.~Jeanjean]{Louis Jeanjean}
	\author[J.~J.~Zhang]{Jianjun Zhang}
	\author[X.~X.~Zhong]{Xuexiu Zhong}

\address[L.~Jeanjean]{\newline\indent Universit\'e Marie et Louis Pasteur
\newline\indent
CNRS, LmB (UMR 6623), F-25000 Besan\c{c}on, France}
\email{\href{mailto:louis.jeanjean@univ-fcomte.fr}{louis.jeanjean@univ-fcomte.fr}}

\address[J.~J.~Zhang]{\newline\indent College of Mathematics and Statistics
\newline\indent
Chongqing Jiaotong University
\newline\indent
Xuefu, Nan'an, 400074, Chongqing, PR China}
\email{\href{mailto:zhangjianjun09@tsinghua.org.cn}{zhangjianjun09@tsinghua.org.cn}}

	\address[X.~X.~Zhong]{\newline\indent South China Research Center for Applied Mathematics and Interdisciplinary Studies \& School of Mathematical Sciences
		\newline\indent
	South China Normal University
		\newline\indent
		Guangzhou 510631, P. R. China}
	 \email{\href{mailto:zhongxuexiu1989@163.com}{zhongxuexiu1989@163.com}}

\thanks{Xuexiu Zhong was supported by the NSFC (No.12271184), Guangdong Basic and Applied Basic Research Foundation (2021A1515010034), Guangzhou Basic and Applied Basic Research Foundation(2024A04J10001). Jianjun Zhang was supported by the NSFC (No.12371109)}

\date{\today}

\begin{document}
\renewcommand{\abstractname}{\normalfont\bfseries Abstract} % ÐÞ¸´ÕªÒª±êÌâ		
\begin{abstract}
This paper establishes optimal existence results and limiting profiles for \textit{energy ground states} of the quasi-linear Schr{\"o}dinger equation
$$
-\Delta u - \Delta(|u|^{2})u + \lambda u = |u|^{p-2}u \quad \text{in } \mathbb{R}^N
$$
with prescribed mass $\int_{\mathbb{R}^N}|u|^2 = a > 0$, in the mass-supercritical case $4 + \frac{4}{N} < p < 2 \cdot 2^*$.

\textbf{Breakthrough in existence theory:} For all dimensions $N \geq 1$, we completely resolve the existence problem: For $1 \leq N \leq 4$, ground states exist for all $a > 0$. For $N \geq 5$, there exists a sharp threshold $a_0 > 0$ such that ground states exist if and only if $a \leq a_0$. This constitutes the optimal existence theory, crucially removing the restrictive condition $p \leq 2^*$ required in all prior works (which limited results to $N \leq 3$).

\textbf{Asymptotic behavior \& new phenomena:} We provide a complete asymptotic analysis of normalized ground states: As $a \to 0^+$, solutions exhibit a novel connection to Serrin-type overdetermined problems. Through a delicate rescaling, profiles converge to the unique positive radial solution of the overdetermined problem (the first such result for quasi-linear equations). As $a \to a^*$ ($a^* = \infty$ for $N \leq 4$; $a^* = a_0$ for $N \geq 5$), solutions converge to distinct limiting profiles depending on dimension and nonlinearity.

Our methods introduce a new constraint approach and unified variational framework for quasi-linear problems with $L^2$-constraints.
\end{abstract}

\maketitle
	
%\tableofcontents

\section{Introduction}
In this paper, we are concerned with quasi-linear Schr\"odinger equations of the form
\begin{equation}\label{eq1.1}
\begin{cases}
i\partial_t\psi+\Delta \psi+\Delta(|\psi|^2)\psi+|\psi|^{p-2}\psi=0, ~\hbox{in}~\R^+\times \R^N,\\
\psi(0,x)=\psi_0(x), ~\hbox{in}~\R^N,
\end{cases}
\end{equation}
where $p\in (2, \frac{4N}{N-2})$ if $N\geq 3$ and $p\in (2,\infty)$ if $N=1,2$, $i$ denotes the imaginary unit and the unknown $\psi:\R^+\times \R^N\rightarrow \mathbb{C}$ is a complex valued function.

Such types of equations appear in various physical fields, for instance in dissipative quantum mechanics, in plasma physics and in fluid mechanics. We refer the readers to \cite{CJS,PSW} and their
 references for more information on the related physical backgrounds.

 From the physical as well as the mathematical point of view, a central issue is the existence and dynamics of standing waves of Eq.\eqref{eq1.1}.
A standing wave of Eq.\eqref{eq1.1} is a solution of the form
$\psi(t,x)=e^{i\lb t}u(x)$, where $\lb\in\R$ is a constant and $u$ is
a time-independent function. Clearly, $\psi$ is a standing wave of
Eq.\eqref{eq1.1} if and only if $u$ is a solution to
\begin{equation}\label{eq1.2}
-\Delta u-\Delta(|u|^2)u+\lambda u =|u|^{p-2}u\ \ \mbox{in} \ \R^N.
%\tag{$A_{\lambda}$}
\end{equation}

Over the past two decades, Eq.\eqref{eq1.2} has received a great deal of attention. Two distinct paths have been explored, one consisting in studying Eq.\eqref{eq1.2} assuming that the frequency $\lb$ is a fixed parameter, the other consisting in searching for solutions to Eq.\eqref{eq1.2} having a prescribed $L^2$ norm.

When $\lambda\in \R$ appears as a fixed parameter in Eq.\eqref{eq1.2}, we refer to \cite{AW2003,CLW2012,CJ,CJS,Iliev1993,LW2003,LW2014,LWW2003,LWW2004,LLW1,LLW2013,LLW2013-2,PSW,RS,ZZ,ZLTZ} and their references therein, for the existence and multiplicity of solutions.
We also refer to \cite{AdaW2012,AW2003,GS2012,Selvitella2011} for uniqueness results.
%of {\it action ground states}.
Formally, a solution to Eq.\eqref{eq1.2} corresponds to a critical point of the {\it action functional}
$$I_\lambda(u):=\frac{1}{2}\|\nabla u\|_2^2+V(u) +\frac{\lambda}{2}\|u\|_2^2-\frac{1}{p}\|u\|_p^p,$$
defined on the natural space
\begin{equation*}
%\beq\lab{eq:Def-X}
X:=\left\{u\in H^1(\R^N)\Big| V(u):=\int_{\R^N}u^2|\nabla u|^2\mathrm{d}x<\infty\right\}.
\end{equation*}
Indeed, it is easy to check that $u$ is a weak solution of Eq.\eqref{eq1.2} if and only if
$$I'_\lambda(u)\phi:=\lim_{t\rightarrow 0^+}\frac{I_\lambda(u+t\phi)-I_\lambda(u)}{t}=0,$$
for every direction $\phi\in C_0^\infty(\R^N,\R)$. However, note that $I_\lambda$ is not differentiable in space $X$ when $N\geq 2$, due to the presence of the term $u\Delta (|u|^2)$, see \cite[page 4]{CaDe}. To overcome this difficulty, various arguments have been developed. First in \cite{LW2003,LWW2004}, solutions to Eq.\eqref{eq1.2} are obtained by minimizing the functional $I_\lambda$ on the constraint
$\displaystyle \left\{u\in X\Big|~~ \|u\|_p^p=1\right\}$.
Then the non-differentiability of $I_\lambda$ essentially does not come into play.

Alternatively, in \cite{CJ, LWW1}, by a change of unknown, Eq.\eqref{eq1.2} is transformed into a semi-linear one, see Subsection \ref{equivalent-problems} for more details. One then speaks of the dual problem.  For this semi-linear equation standard variational methods can be applied to yield a solution. Also in \cite{LWW2} the authors have developed an approach which works for more general quasi-linear equations and reduces the search of solutions to Eq.\eqref{eq1.2} to the problem of showing that
$I_{\lambda}$ has a global minimizer on a Nehari manifold.
Since this pioneering work, the study of Eq.\eqref{eq1.2}  has been the subject of an abundant literature that addresses questions of multiplicity, concentration and critical exponent.
In particular, starting from \cite{LLW2013}, in a series of papers \cite{LW2014,LLW1,LLW2013-2,LLW2} Liu, Liu and Wang have developed a perturbation method for studying the existence and multiplicity of solutions for a general class of quasi-linear elliptic equations including the above discussed model equations. Finally, we mention the recent work \cite{Genoud-Nodari-2024} where it is proved, that for each $\lambda >0$, there exists a unique positive solution to Eq.\eqref{eq1.2}. In addition, it is non-degenerate, and lie on a $C^1$ branch.

In this paper, we consider the other approach and search for solutions $u\in X$ to Eq.\eqref{eq1.2} with a prescribed mass. Namely, such that
\beq\lab{eq:contraint-1}
u\in \mS_a:=\left\{u\in X\Big|~~ \|u\|_2^2=a\right\}
\eeq
where $a>0$ is given.
Then the frequency $\lambda$ will appear as a Lagrange multiplier. This kind of solution is called a normalized solution in the literature and the search for it is motivated by the fact that the $L^2$ norm has, depending on the underlying model, a strong physical meaning.
Solutions to Eq.\eqref{eq1.2}-Eq.\eqref{eq:contraint-1} can be obtained as critical points
of the {\it energy functional}
\begin{equation*}
%\lab{eq:functional}
I(u)=\frac{1}{2}\|\nabla u\|_2^2+V(u)-\frac{1}{p}\|u\|_p^p
\end{equation*}
constrained on $\mS_a$.
Due to a Gagliardo-Nirenberg
type inequality, see Eq.\eqref{eq2.1}, the
functional $I$ is bounded from below on $\mS_a$ if $p<4+\frac{4}{N}$ and
is unbounded from below on $\mS_a$ if $p>4+\frac{4}{N}$. The threshold
number $4+\frac{4}{N}$ is called the mass critical exponent for
Eq.\eqref{eq1.2} in the literature.

Along this approach a natural notion is the one of \textit{least energy solution}. A \textit{least energy solution} is a solution $u \in \mathcal{S}_a$ that minimizes the {\it energy functional} $I$ among all solutions to Eq.\eqref{eq1.2} subject to the mass constraint Eq.\eqref{eq:contraint-1}.
This definition remains meaningful even when $I$ is unbounded below on $\mathcal{S}_a$. Moreover, when $I$ is bounded below on $\mathcal{S}_a$, it makes sense even if the $\inf_{u \in \mathcal{S}_a} I(u)$ is not attained.

%The natural counterpart to \textit{action ground states} in this constrained setting is the concept of \textit{energy ground states}. An \textit{energy ground state} is a solution $u \in \mathcal{S}_a$ that minimizes the energy functional $I$ among all solutions to \eqref{eq1.2} subject to the mass constraint \eqref{eq:contraint-1}.
%This definition remains meaningful even when $I$ is unbounded below on $\mathcal{S}_a$. Moreover, when $I$ is bounded below on $\mathcal{S}_a$, it is also make sense even if the infimum $\inf_{u \in \mathcal{S}_a} I(u)$ is not attained.

%The equivalent of the {\it action ground states} are now the {\it energy ground states}. An {\it energy ground state} is a $u \in \mS_a$ which minimizes among all the solutions to Eq.\eqref{eq1.2}-Eq.\eqref{eq:contraint-1} the energy functional $I$.
%Note that this definition, makes sense even if $I$ is unbounded from below on $S_a$, and also that when $I$ is bounded from below on $S_a$ an {\it energy ground state} can exist even if the infimum of $I$ on $S_a$ is not achieved.

The mass subcritical case (i.e., $2<p<4+\frac{4}{N}$) has now been widely studied. In \cite{CJS,JLuo} {\it least energy solutions} are obtained by minimizing the functional $I$ on $\mS_a$, see \cite[Lemma 1.1]{JLWang} for results in that direction. Critical points which are not global minimum have also been obtained.
In \cite{JLWang}, for $2+\frac{4}{N}<p<4+\frac{4}{N}$, a perturbation method is applied to study the existence and multiplicity for suitable range of $a$.
In \cite{ZLW}, the authors adopted the dual approach to construct infinitely many pairs of normalized solutions when $2<p<2+\frac{4}{N}$. Also, when $2+\frac{4}{N}\leq p<4+\frac{4}{N}$, $k$ pairs of normalized solutions were obtained in \cite{ZLW}, for any $k \in \N$, assuming that $a >a_k$ for some $a_k$.

In the mass critical case (i.e., $p=4+\frac{4}{N}$), the infimum of $I$ on $\mS_a$ is never achieved but, for  $a>0$ sufficiently large, one may still find a {\it least energy solution} \cite{ZCW}. Also, still for $a>0$ sufficiently large, the existence (of mountain pass solution) and concentration behavior of normalized solutions of  Eq.\eqref{eq1.2}-Eq.\eqref{eq:contraint-1} are established in \cite{LZou}.

 The mass supercritical
case (i.e., $4+\frac{4}{N}<p<\frac{4N}{N-2}$), is still essentially untouched. In the works
\cite{LZou,ML}, existence and multiplicity results for  Eq.\eqref{eq1.2}-Eq.\eqref{eq:contraint-1} are obtained using a  perturbation method.  In \cite{ZCW}, the authors prove the existence of {\it least energy solutions} by minimizing the functional $I\big|_{\mS_a}$ on the so-called Pohozaev manifold
\begin{equation*}
\mathcal{P}:=\{u\in X\backslash\{0\}: P(u)=0\},
\end{equation*}
where $$P(u):=\|\nabla u\|_2^2+(N+2)V(u)-\frac{(p-2)N}{2p}\|u\|_p^p.$$
As we shall see, $\mathcal{P}$ is a natural constraint which contains all the critical points of $I$ constrained to $\mS_a$.

%In the approach of \cite{ZCW}, the non-differentiability of $I$ essentially does not come into play.
%For the \textcolor{red}{mass-}mixed case, we also note that replacing $|u|^{p-2}u$ by a double power $f(u)=\tau |u|^{q-2}u+|u|^{p-2}u$ with $2<q<2+\frac{4}{N}, 4+\frac{4}{N}<p\begin{cases}<\infty, N=1,2,\\ \leq 2^*, N=3 \end{cases}$, the authors in \cite{ML} also adopted a perturbation argument to study the existence of normalized ground state solution and the mountain pass solution.

However, all previous studies \cite{LZou,ZCW,ML} require the restrictive assumption $p \leq 2^*$ (or $p < 2^*$). This limitation has critical consequences:
\begin{itemize}
    \item Results are confined to dimensions $N \leq 3$.
    \item Even for $N=3$, the full exponent range $4 + \frac{4}{N} < p < 2\cdot 2^* = \frac{4N}{N-2}$ remains inaccessible.
\end{itemize}

The fundamental barrier stems from the necessity in \cite{LZou,ZCW,ML} to ensure the strict positivity of the Lagrange multiplier $\lambda$ (a crucial element for controlling potential compactness loss). Specifically, their analysis relies exclusively on two identities for solutions $u$:
\begin{equation}\label{eq:identities}
\begin{cases}
\text{Nehari identity:} & \|\nabla u\|_2^2 + 4V(u) + \lambda\|u\|_2^2 = \|u\|_p^p,\\
\text{Pohozaev identity:} & P(u) = 0.
\end{cases}
\end{equation}
Combining these yields the expression for $\lambda$:
\begin{equation}\label{eq:20241025-1321}
\begin{split}
\lambda \|u\|_2^2 &=\|u\|_p^p-\left[\|\nabla u\|_2^2 + 4V(u)\right]\\
&= \frac{2p}{(p-2)N} \left[\|\nabla u\|_2^2 + (N+2)V(u)\right] - \left[\|\nabla u\|_2^2 + 4V(u)\right] \\
&= \underbrace{\frac{2p - (p-2)N}{(p-2)N}}_{\text{require $\geq 0$}} \|\nabla u\|_2^2 + \frac{2p(N+2) - 4(p-2)N}{(p-2)N} V(u).
\end{split}
\end{equation}
The non-negativity requirement for the coefficient $\frac{2p-(p-2)N}{(p-2)N}$ forces $p \leq 2^*$.

\textbf{Our fundamental breakthrough:} Focusing on the existence of {\it energy ground states}, namely {\it least energy solutions} which minimizes
$I\big|_{\mS_a}$ on $\mathcal{P}$,
this work eliminates this long-standing restriction on $p$. We establish comprehensive results for:
\begin{itemize}
    \item \textit{All dimensions} $N \geq 1$.
    \item The \textit{complete range} $p \in \left(4 + \frac{4}{N}, 2\cdot 2^*\right)$.
\end{itemize}
Additionally, we conduct a thorough asymptotic analysis of {\it energy ground states} as their $L^2$-norm approaches critical values ($a \to 0^+$ and $a \to a^*$ where $a^* = \infty$ for $N \leq 4$, $a^* = a_0$ for $N \geq 5$).

\section{Main results}
\subsection{Existence of energy ground states}
We will work on the space $X$.
Note that $X$ is not a vector space since it is not closed under the
sum. Nevertheless, $X$ is a complete metric space with the distance
$\mathrm{d}_X(u, v)=\|u-v\|_{H^1(\R^N)}+\|\nabla (u^2)-\nabla (v^2)\|_{L^2(\R^N)}$.
Due to the presence of $V(u)$, $I$ is not differentiable in $X$ when $N\geq 2$, and powerful variational techniques, such as those involving Palais-Smale sequences, cannot be used.
To remedy this, inspired by \cite{ZCW}, we develop a minimization approach that uses only elementary arguments and prove, a posteriori, that a minimum of $I$ on $\mathcal{P} \cap \mS_a$ is indeed a critical point of $I$ restricted to $\mS_a$. However, the strategy of the proofs in \cite{ZCW} and in the present work is essentially distinct.

Instead of working directly on the constraint $\mathcal{P}\cap\mS_a$, as in \cite{ZCW}, inspired by the work \cite{Jeanjean-Zhang-Zhong-2024}, see also \cite{Clarke1976,DHZ-2026,BM} for the introduction of this strategy, we first consider a relaxed problem. For any $a>0$, we introduce the set
$D_a:=\{u\in X: \|u\|_2^2\leq a\}$
and we define $\mathcal{P}_a:=\mathcal{P}\cap D_a.$
The set $\mathcal{P}_a$ is a natural constraint since, as a consequence of a Pohozaev identity, any non-trivial critical point of $I$ constrained to $D_a$ belongs to $\mathcal{P}_a$.

The relaxed problem we will consider is:
\begin{align}
\lab{eq:20220930-e1}
\mbox{to find a minimizer}~u\in \mathcal{P}_a\mbox{ of } M_a \mbox{ where } M_{a}:=\inf_{\mathcal{P}_{a}}I(u).
\end{align}
By Lemma \ref{lemma:20240830-1318}, a minimizer of $M_a$ will prove to be a critical point of $I_{|_{D_a}}$ and, by the above consideration, it will minimize $I$ among all non-trivial critical points of $I_{|_{D_a}}$.

Once this problem is settled, we will show that $u\in \mS_a$, and thus that it is an {\it energy ground state} for the initial problem.

Here comes our first main result.
\bt\lab{main-th-1}
Let $1\leq N\leq 4$ and assume that $p\in (4+\frac{4}{N}, 2\cdot 2^*)$, then the following conclusions hold true.
\begin{itemize}
\item[(i)] $M_a$ is continuous in $(0, +\infty)$;
\item[(ii)] $M_a$ decreases strictly in $(0, +\infty)$;
\item[(iii)] Eq.\eqref{eq1.2}-Eq.\eqref{eq:contraint-1} has an energy ground state and, up to translation and sign change, any energy ground state is positive, radial, and strictly decreasing  with respect to $r = |x|$. In addition, the associated Lagrange parameter is positive.
\end{itemize}
\et

For $1\leq N\leq 4$, by Theorem \ref{main-th-1} above, the existence of an {\it energy ground state} holds for all prescribed mass $a>0$. The situation when $N\geq 5$ is different. In order to present it, we first recall the following result that can be found  in \cite{Genoud-Nodari-2024}.

\bo\lab{remark:20241001-1250}
For $N\geq 3$ with $p\in (2^*, 2\cdot 2^*)$, there exists a unique positive, radial solution $u\in D_{0}^{1,2}(\R^N) \cap L^p(\R^N)$ of the equation
    \beq\lab{eq:20240830-1747}
-\Delta u-\Delta(|u|^2)u  =|u|^{p-2}u\ \ \mbox{in} \ \R^N.
\eeq
Furthermore, $u(x)=O(|x|^{-(N-2)})$ as $|x|\rightarrow \infty$. In particular, $u\in L^2(\R^N)$ when $N\geq 5$.
\eo
Note that Eq.\eqref{eq:20240830-1747} corresponds to the so-called {\it zero-mass case} in Eq.\eqref{eq1.2}. As indicated in \cite{Genoud-Nodari-2024}, the existence part follows directly from the results in  \cite[Section 5]{BL} and the uniqueness part from \cite[Theorem 2]{Tang-2001}.

%Let $u=\varphi(v)$ be the inverse function of
%\beq\lab{eq:20241109-2106}
%v=\int_0^u\sqrt{1+2t^2}\,dt
%=\frac12 u\sqrt{1+2u^2}+\frac{\sqrt{2}}{4}\ln\big(\sqrt{2}u+\sqrt{1+2u^2}\big).
%\eeq
%Then one can also use the classical result in \cite{SeT} to obtain the uniqueness of $v$ for the semilinear equation
%\beq\lab{eq:20241001-1300}
%-\Delta v=\varphi(v)^{p-1}\varphi'(v)~\hbox{in}~\R^N,
%\eeq
% which equivalents to the uniqueness of $u$.

Let $N\geq 5$ and $4+\frac{4}{N}<p<2\cdot 2^*$, then it holds that $p\in (2^*, 2\cdot 2^*)$. Let $u_0\in D_{0}^{1,2}(\R^N) \cap L^p(\R^N)$ be given by Proposition \ref{remark:20241001-1250}, and define
\beq\lab{eq:20241001-1308}
a_0:=\|u_0\|_2^2.
\eeq
When $N \geq 5$, our result regarding the existence of {\it energy ground states} is as follows.
\bt\lab{main-th-2}
Assume that $N\geq 5$ and $p\in (4+\frac{4}{N}, 2\cdot 2^*)$. Let $a_0$ be given by Eq.\eqref{eq:20241001-1308}, then the following conclusions hold true.
\begin{itemize}
\item[(i)] $M_a$ is continuous in $(0, +\infty)$;
\item[(ii)] $M_a$ decreases strictly in $(0, a_0]$ and $M_a\equiv M_{a_0}$ in $[a_0, +\infty)$;
\item[(iii)] If  $a\in (0, a_0)$, Eq.\eqref{eq1.2}-Eq.\eqref{eq:contraint-1} has an {\it energy ground state} and, up to translation and sign change, any energy ground state is positive, radial, and strictly decreasing  with respect to $r = |x|$. The associated Lagrange parameter is positive;
\item[(iv)] If $a =a_0$, Eq.\eqref{eq1.2}-Eq.\eqref{eq:contraint-1} admits $u_0$ as an {\it energy ground state} and it is, up to translation and sign change, the unique energy ground state. The associated Lagrange parameter is zero;
\item[(v)] If $a > a_0$, Eq.\eqref{eq1.2}-Eq.\eqref{eq:contraint-1} has no {\it energy ground state}.
\end{itemize}
\et

\br\lab{remark:openproble}
When there is no energy ground state, one cannot conclude that {\it least energy solutions} do not exist, but only that, if they do exist, they lay at a higher energy level. Note also that a {\it least energy solution} is not necessarily a positive function, see \cite[Remark 1.6]{DDGS} in that direction. The existence, or non-existence, of
{\it least energy solutions} whenever {\it energy ground states} do not exist appears fully open.
\er

Our existence results for {\it energy ground states} are summarized as follows.
\begin{table}[H]
\centering
%\caption{Existence and nonexistence of minimizer of $M_a$}
\medskip
\begin{tabular}{|c|c|c|c|}
\hline
\text{Dimension }&\text{Range of mass}&\text{Value of $\lambda$}&\text{Existence}\\
\hline
\text{$1\leq N\leq 4$}&\text{$a\in (0,+\infty)$}&\text{$\lambda_a>0$}&\text{Yes}\\
\hline
\multirow{3}{2cm}{$N\geq 5$}
&\text{$a\in (0,a_0)$}&\text{$\lambda_a>0$}&\text{Yes}\\
\cline{2-4}
&\text{$a=a_0$}&\text{$\lambda_a=0$}&\text{Yes}\\
\cline{2-4}
&\text{$a\in(a_0,+\infty)$}&/&\text{No}\\
\hline
\multicolumn{4}{|c|}{$N\geq 1, 4+\frac{4}{N}<p<2\cdot 2^*$}\\
\hline
\end{tabular}
\end{table}

\br\lab{remark:threshold}
Comparing Theorem \ref{main-th-1} and
Theorem \ref{main-th-2}, one finds that the existence of a positive (radial) solution in $L^2(\R^N)$ to Eq.\eqref{eq:20240830-1747} crucially influences the existence of {\it energy ground states}. This role played by a solution of the corresponding {\it zero-mass} equation has already been observed in \cite{BFG} on a Schr\"odinger equation with competing nonlinearities. We also refer to \cite{BJL} for partial results in that direction.
\er

\br\lab{remark:stability}
Having proved the existence of {\it energy ground states}, a natural question is the orbital stability, or instability, of the associated standing waves in Eq.\eqref{eq1.1}. In view of the classical results of \cite{BeCa}, sees also \cite{BJL,Soave}, the variational characterization of the energy ground states in Theorems \ref{main-th-1} and \ref{main-th-2}, as minimizer on the constraint $\mathcal{P} \cap \mS_a$, is a strong indication that these solutions are unstable by blow-up. However, we cannot directly conclude this since a local well-posedness theory of Eq.\eqref{eq1.1}  in the energy space $X$ is still to be established.  See \cite[Theorem 1.1]{CJS} for a partial result in that direction.
\er

\br\lab{remark:action versus energy ground state}
Among the solutions to Eq.\eqref{eq1.2} when $\lambda \in \R$ is fixed, particular ones correspond to {\it action ground states}, namely solutions to Eq.\eqref{eq1.2} which minimize $I_{\lambda}$ among all non-trivial solutions. In general, in related settings, it is a difficult issue to determine a link between {\it action ground states} and {\it least energy solutions}, only partial results are known \cite{DDGS, DoSeTi, JeLu}. Here, thanks to the uniqueness of positive solutions when $\lambda >0$ is fixed obtained in \cite{Genoud-Nodari-2024} (leading, in particular, to the fact that they are {\it action ground states}), we  conclude that our {\it energy ground states} are also {\it action ground states}.
\er

\subsection{Limiting profile of the energy ground states}

In \cite{Genoud-Nodari-2024}, the question of the orbital stability of the {\it action ground states} is addressed, in the spirit of the work of Grillakis-Shatah-Strauss
\cite{GrShSt1}. This is possible because, for any $\lambda \in (0, \infty)$,  the unique {\it action ground state} $u_{\lambda}$ lies on a $C^1$ curve. From \cite{GrShSt1} it is known that the standing wave is expected to be orbitally stable if $\lambda \mapsto ||u_{\lambda}||_2$ is strictly increasing and unstable when $\lambda \mapsto ||u_{\lambda}||_2$ is strictly decreasing. Because the full problem of stability, namely of $u_{\lambda}$ for any $\lambda \in (0, \infty)$, seems out of reach so far, one may consider what happens in some limiting regimes.  In \cite{Genoud-Nodari-2024} a very complete description of the behavior of $\lambda \mapsto ||u_{\lambda}||_2$ is given in the limit $\lambda \to 0^+$. To study the behavior of this function a key step is to identify a limit equation, having a unique solution, to which $u_{\lambda}$, or some rescaling of it, converge.

Motivated, in particular, by the work \cite{Genoud-Nodari-2024}, we now consider what happens to our {\it energy ground states} in the limit $a \to 0^+$ and $a \to  a^*$ respectively, where
\beq\lab{eq:20241002-2126}
a^*:=\begin{cases}
+\infty,\quad&\hbox{if}~1\leq N\leq 4,\\
a_0,\quad&\hbox{if}~N\geq 5.
\end{cases}
\eeq

First we consider the case where the mass $a$ goes to zero.
%In that direction we have the following result.
\bt\lab{th:20240905-2122}
Assume that $N\geq 1,$ $p\in (4+\frac{4}{N}, 2\cdot 2^*)$.
Let $\{u_n\}_{n=1}^{\infty}$, $u_n\in \mS_{a_n}$, be a sequence of {\it energy ground states} to Eq.\eqref{eq1.2}-Eq.\eqref{eq:contraint-1} with $a_n\rightarrow 0^+$. Then the associated Lagrange multiplier satisfies $\lambda_n\rightarrow + \infty $ as $n\rightarrow \infty$. Furthermore,  defining
\begin{equation*}
%\beq\lab{eq:20241109-2106}
v_n=\int_0^{u_n}\sqrt{1+2t^2}\,\mathrm{d}t
=\frac12 u_n\sqrt{1+2u_n^2}+\frac{\sqrt{2}}{4}\ln\big(\sqrt{2}u_n+\sqrt{1+2u_n^2}\big)
%\eeq
\end{equation*}
and
$\displaystyle \tilde{v}_n(x):=\lambda_{n}^{-\frac{2}{p-2}}v_n\left(\frac{x}{\lambda_{n}^{\frac{p-4}{2(p-2)}}}\right),$
it holds that,  up to a subsequence, $\tilde{v}_n\rightarrow \tilde{v}$ in $C^{1,\alpha}(\R^N)\cap C_{\text{loc}}^{2}(\Omega)$ with $\Omega:=\text{supp}~\tilde{v}$, where $\tilde{v}$ is  the unique positive radial solution to the following overdetermined problem
\beq\lab{eq:20240905-1912}
\begin{cases}
-\Delta \tilde{v}=-\frac{\sqrt{2}}{2}+2^{\frac{p-4}{4}}\tilde{v}^{\frac{p-2}{2}},\\
\tilde{v}>0~\hbox{in}~\Omega,~~\tilde{v}=\frac{\partial \tilde{v}}{\partial \nu}=0~\hbox{on}~\partial \Omega,\\
\hbox{where $\Omega=B_R$ is an open ball in $\R^N$ with radius $R>0$}.
\end{cases}
\eeq
Precisely,  there exists some $R>0$ such that $\tilde{v}$ is supported in $B_R$ and it satisfies the above Dirichlet-Neumann boundary problem.
\et

The uniqueness of the positive radial solutions to the overdetermined Problem \eqref{eq:20240905-1912} can be deduced from existing results. When $N\geq 3$, Problem \eqref{eq:20240905-1912} is a special case of the general result presented in \cite[Theorem 2]{SeT}. When $N=1,2$ we make use of the very recent result \cite[Theorem 1.4]{Jeanjean-Zhang-Zhong-2025}.

%When $N=2$, we shall make use of the result \cite[Theorem 3]{Pucci-Serrin1998}. The case where $N=1$ was apparently still open, but it is now treated in the work \cite{Jeanjean-Zhang-Zhong-2025}.

Regarding the convergence directly on $u_n$-functions, we establish the following result.
\bt\lab{th:20250101-1347}
Assume that $N\geq 1,$ $p\in (4+\frac{4}{N}, 2\cdot 2^*)$. As $a_n\rightarrow 0^+$, let $\{u_n\}, \{\lambda_n\}, \tilde{v}, \Omega$ be given in Theorem \ref{th:20240905-2122}.
Define
$\displaystyle
\bar{u}_n(x):=\lambda_{n}^{-\frac{1}{p-2}}u_n\left(\frac{x}{\lambda_{n}^{\frac{p-4}{2(p-2)}}}\right).
$
Then $\bar{u}_n\rightarrow \bar{u}$ in $C_{\text{loc}}^{2}(\Omega)$ and $\bar{u}$ is the unique positive radial solution to
\beq\lab{eq:20250101-1346}
\begin{cases}
- \Delta(|\bar{u}|^2)+1=\bar{u}^{p-2}~\hbox{in}~\Omega,\\
\bar{u}(0)=2^{\frac{1}{4}} \sqrt{\tilde{v}(0)}, \, \bar{u}'(0)=0.
\end{cases}
\eeq
In particular, we have that $\bar{u}=2^{\frac{1}{4}}\sqrt{\tilde{v}}$.
\et

%\medskip

We now consider the limit $a \uparrow a^*$.
\bt\lab{th:20240928-2149}
Let $\{u_n\}_{n=1}^{\infty}$, $u_n\in \mS_{a_n}$, be a sequence of {\it energy ground states} to Eq.\eqref{eq1.2}-Eq.\eqref{eq:contraint-1} with $a_n \uparrow a^*$. Then the associated Lagrange multiplier satisfies $\lambda_n\rightarrow 0^+$ as $n\rightarrow \infty$. Furthermore, the following conclusions hold true.
\begin{itemize}
\item[(i)] If $N=2$ or $N=3$, with $4+\frac{4}{N}<p<2^*$ then the rescaled functions
$\displaystyle
\lambda_{n}^{-\frac{1}{p-2}}u_n\left(\frac{x}{\sqrt{\lambda_{n}}}\right),
$
converge in $H^1(\R^N) \cap C^2(\R^N)$ to the unique positive radial solution to
\begin{equation*}
-\Delta W+W=W^{p-1}~\hbox{in}~\R^N, ~\lim_{|x|\rightarrow \infty}W(x)=0.
\end{equation*}
\item[(ii)]If $N=3$ and $p=2^*$, there exists $\mu_n
=O(\lambda_{n}^{-\frac{1}{4}})\rightarrow +\infty$ such that the rescaled functions $\displaystyle \mu_{n}^{\frac{1}{2}} u_n(\mu_n x)$ converge
in $D_{0}^{1,2}(\R^3)\cap C^2(\R^3)$ to the Talenti bubble (up to a multiplier) $U(x)=3^{\frac{1}{4}}\big(1+|x|^2\big)^{-\frac{1}{2}},$
which is the unique (up to dilatation) positive radial decreasing solution to
\begin{equation*}
-\Delta U=U^{5}~\hbox{in}~\R^3, ~\lim_{|x|\rightarrow \infty}U(x)=0.
\end{equation*}
\item[(iii)] If $N=3$ with $2^* < p < 2\cdot 2^*$ or $N\geq 4$ with  $ 4+ \frac{4}{N}< p < 2\cdot 2^*$ , as $a_n \uparrow a^*$, $u_n$ converges in $D_{0}^{1,2}(\R^N)\cap C^2(\R^N)$ to the unique positive radial decreasing  solution $u_0\in D_{0}^{1,2}(\R^N)\cap L^{p}(\R^N)$ of  Eq.\eqref{eq:20240830-1747}. In particular, since $u_0=O(|x|^{-(N-2)})$ as $|x|\rightarrow \infty$, we have that $u_n \to u_0$ in $L^2(\R^N)$ as $a_n \uparrow a^*$ when $N \geq 5$.
\end{itemize}
\et

As we shall see, Theorem \ref{th:20240928-2149} follows directly from \cite[Theorem 1.3]{Genoud-Nodari-2024} once it is proven that $\lambda_n \to 0^+$ as $a_n \uparrow a^*$.

Note that the case $N=1$ is not treated in Theorem \ref{th:20240928-2149} as it is not in \cite{Genoud-Nodari-2024}. When $N=1$, the work \cite{Iliev1993} already provides comprehensive results concerning the question of stability. Actually, using arguments from \cite{Jeanjean2024} and from the proof of Theorem \ref{th:20240905-2122}, it can be shown that Theorem \ref{th:20240928-2149} (i) also holds when $N=1$.
We do not include a proof of this result to keep the paper within reasonable length.

%%%%%%%%%%%%%%%%%%%%%%%%%%%%%%%%%%%%%%%%%%%%%%%%%%%%%%%%%%%%%%%%%%%%%%%%%%%
%%%%%%%%%%%%%%%%%%%%%%%%%%%%%%%%%%%%%%%%%%%%%%%%%%%%%%%%%%%%%%%%%%%%%%%%%%%

\section{Preliminaries}

\subsection{A Gagliardo-Nirenberg type inequality}

We recall the following Gagliardo-Nirenberg type inequality: There is some positive constant $C(N,p)$ such that
%\begin{equation}\label{eq2.1}
%\|u\|_p^p\leq C(N,p)\|u\|_{2}^{\frac{4N-p(N-2)}{N+2}}
%\left(\int_{\R^N} |u|^2\, |\nabla u|^2{\rm d}x\right)^{\frac{N(p-2)}{2(N+2)}}, \forall u\in X, p\in (2, 2\cdot2^*).
%\end{equation}
\begin{equation}\label{eq2.1}
\int_{\R^N}|u|^p\,{\rm d}x
\leq C(N, p) \left(\int_{\R^N} u^2\,{\rm d}x\right)^{\frac{4N-p(N-2)}{2(N+2)}}
\left(\int_{\R^N} |u|^2\, |\nabla u|^2{\rm d}x\right)^{\frac{N(p-2)}{2(N+2)}}, 
\end{equation}
for $u \in X$ and $p \in (2, 2\cdot2^*)$.

For a proof of Eq.\eqref{eq2.1} we refer to \cite[Lemma 4.2]{CJS} or to \cite[Lemma 4.14]{Genoud-Nodari-2024}.

\subsection{Two equivalent problems}\label{equivalent-problems}

It has been observed
%for a long time
that the search of non-negative solutions to Eq.\eqref{eq1.2} is equivalent to the search of non-negative solutions to
\beq\lab{eq:20240904-0823}
-\Delta v=-\lambda \varphi(v)\varphi'(v)+\varphi(v)^{p-1}\varphi'(v)~\hbox{in}~\R^N,
\eeq
where $\varphi$ is the unique solution of the Cauchy problem
$$
\varphi'(s)=\frac{1}{\sqrt{1+2\varphi^2(s)}},\quad\varphi(0)=0
$$
on $[0, + \infty)$ and $\varphi(t) = - \varphi(-t)$ on $(- \infty, 0].$
Indeed, by \cite{CJ} it is known that any non-negative solution $v\in H^1(\R^N)$ of Eq.\eqref{eq:20240904-0823} gives rise to a non-negative solution $u:=\varphi(v)$ of Eq.\eqref{eq1.2}. Reciprocally, from \cite[Lemmas 2.6, 2.8, 2.10]{AdaW2012} we know that for any non-negative solution $u$ of Eq.\eqref{eq1.2}, $v:=\varphi^{-1}(u)$ is a non-negative solution of Eq.\eqref{eq:20240904-0823}.
%Clearly also a non-negative solution of  Eq.\eqref{eq1.2} is transformed into a non-negative solution of Eq.\eqref{eq:20240904-0823}.

For future reference we recall some properties of $\varphi$. First note that
\begin{equation}\label{eq3.8}
\varphi'(0)=\lim_{t\to 0}\frac{\varphi(t)}{t}=1,\ \
\lim_{t\to+\infty}\frac{\varphi(t)}{\sqrt{t}}=2^{\frac14},
\end{equation}
and for $t\in\R$
\begin{equation}\label{eq3.9}
0<\varphi'(t)\leq 1,\ \
|\varphi(t)|\leq\min\big\{|t|,\,2^{\frac14}\sqrt{|t|}\big\},\ \
\frac{1}{2} \varphi^2(t)\leq \varphi(t)\varphi'(t)t\leq \varphi^2(t).
\end{equation}
See \cite[Lemma 1.1 and Lemma 1.2]{CJ} for a proof of Eq.\eqref{eq3.8} and  Eq.\eqref{eq3.9}.

Since $u=\varphi(v)$, clearly $u\rightarrow +\infty$ if and only if $v\rightarrow +\infty$. So, for $u\rightarrow +\infty$ (or saying $v\rightarrow +\infty$), we have that
$\sqrt{1+2u^2}\sim \sqrt{2}u =\sqrt{2}\varphi(v)\sim 2^{\frac{3}{4}}v^{\frac{1}{2}}.$
Now, since
$1=\varphi'(v)\frac{\mathrm{d}v}{\mathrm{d}u}=\varphi'(v)\sqrt{1+2u^2},$
%$$1=\varphi'(v)\frac{\mathrm{d}v}{\mathrm{d}u}=\varphi'(v)\sqrt{1+2u^2},$$
it follows that
\beq\lab{eq:20240904-1055}
\lim_{s\rightarrow +\infty}\varphi'(s)s^{\frac{1}{2}}=2^{-\frac{3}{4}}.
\eeq
\begin{lemma}\label{lemma:20260413-2245}
Let $v\in H^1(\mathbb{R}^N)\backslash \{0\}$ be a non-negative solution to Eq.\eqref{eq:20240904-0823}. Then $v \in C^2(\mathbb{R}^N)$ and $v>0$. In addition,
if $\lambda>0$, then up to a translation, $v$ is radial and strictly decreasing w.r.t. $r=|x|$.
\end{lemma}
\begin{proof}
By Moser iteration, $v \in L^\infty(\mathbb{R}^N)$, so the right-hand side of Eq.\eqref{eq:20240904-0823} is in $L^\infty(\mathbb{R}^N)$. Standard elliptic regularity yields $v \in C^2(\mathbb{R}^N)$, and the strong maximum principle implies $v > 0$ in $\mathbb{R}^N$.
Now define $f(s):=-\lambda \varphi(s)+\varphi(s)^{p-1}\varphi'(s)$. Assuming $\lambda>0$, one can check that $f'(s)=-\lambda(1+o(s))<0$
for sufficiently small $s>0$. So, if $N\geq 2$, by \cite[Theorem 1]{LN-1993}, up to a translation, all positive classical solutions to Eq.\eqref{eq:20240904-0823} are  radial and decrease strictly with respect to $r=|x|$.
If $N=1$, by $p>4+\frac{4}{N}$, it is easy to see that
\begin{equation}\label{eq:20260413-2253}
\xi_0:=\inf\{\xi>0; F(\xi)=0\}~\hbox{exists, and}~\xi_0>0.
\end{equation}
By the necessary condition in \cite[Theorem 5]{BL} (due to the existence of $v$), we have that $f(\xi_0)>0$. Then by \cite[Theorem 5]{BL} again, Eq.\eqref{eq:20240904-0823} has, up to translations, a unique positive solution. This solution is radial with $u(0)=\xi_0$ and $u'(x)<0, x>0$.
\end{proof}

%\br\lab{remark:20250121-1047}
%Since the search for solutions $u\in X$ to Eq.\eqref{eq1.2} is equivalent to the one of solutions $v\in H^1(\R^N)$ to Eq.\eqref{eq:20240904-0823}, the techniques employed on the fixed frequency problem are quite similar to those in the semilinear case. However, when the focus shifts to the existence of normalized solutions, the situation changes significantly. The previous $L^2$ constraint transforms into a complex constraint. Consequently, even when it is rephrased in the frame the Eq.\eqref{eq:20240904-0823}, the relevant Gagliardo-Nirenberg inequality cannot be effectively applied under this constraint. This is in particular why we do not use this transformation to study the existence of normalized solutions. Nevertheless, we sometimes rely on this equivalent problem to analyze and derive certain properties of the related solutions.
%\er

\subsection{Uniqueness of positive radial solutions to a class of overdetermined problems}\label{uniqueness}

For the uniqueness of positive radial solutions to Problem \eqref{eq:20240905-1912}, in the case $N\geq 3$, we make use of \cite[Theorem 2]{SeT} and for $N=1,2$ of \cite[Theorem 1.4]{Jeanjean-Zhang-Zhong-2025}. \\ Indeed, when $1 < m <N$, under the assumptions: there exists some $b>0$ such that
\begin{itemize}
\item[(H1)] $f$ is continuous on $(0,\infty)$, with $f(s)\leq 0$ on $(0,b]$ and $f(s)>0$ for $s>b$;
\item[(H2)] $f\in C^1(b,\infty)$, with $g(s)=sf'(s)/ f(s)$ non-increasing on $(b,\infty)$,
\end{itemize}
Serrin and Tang proved the uniqueness of radial solutions of the homogeneous Dirichlet-Neumann free boundary problem
\beq\lab{eq:20240923-0834}
\begin{cases}
-\Delta_m u=f(u), N>m>1,\\
u>0~\hbox{in}~B_R, u=\frac{\partial u}{\partial\nu}=0~\hbox{on}~\partial B_R,
\end{cases}
\eeq
where $B_R$ is an open ball in $\R^N$ with radius $R>0$. However, the low dimension case $1 \leq N \leq m, m >1$ is left as an open problem (see \cite[Subsection 6.2]{SeT}). Recently, it has been resolved by Pucci, Zhang and Zhong. In \cite{Jeanjean-Zhang-Zhong-2025} the following result is established.

\bt\label{th:20240923-1010}
Let $m>1,$ $1 \leq N \leq m$ and $(H1)$ hold.
\begin{itemize}
\item[(i)] If $N=1$ and Eq.\eqref{eq:20240923-0834} admits a radial solution $u$, then $u$ is unique and $u(0) >b$ is uniquely determined by $F(u(0)) =0$, where $F(t) = \int_{0}^{s}f(s) \mathrm{d}s$.
\item[(ii)] Assume further $(H2)$ and \\
$(H3)$ $g(u) > -1$ for $u >b$, \\
when $N \geq 2$, then Eq.\eqref{eq:20240923-0834} admits at most one radial solution.
\end{itemize}
\et

Let us show that the results of \cite{SeT,Jeanjean-Zhang-Zhong-2025} apply to Problem \eqref{eq:20240905-1912}.   In this case
$\displaystyle f(s):=-\frac{\sqrt{2}}{2}+2^{\frac{p-4}{4}} s^{\frac{p-2}{2}}$ and thus the following properties hold:
\begin{itemize}
\item[(i)] $f$ is continuous on $(0,\infty)$, with $f(s)\leq 0$ on $(0,\frac{\sqrt{2}}{2}]$ and $f(s)>0$ for $s>\frac{\sqrt{2}}{2}$;
\item[(ii)]$f\in C^1(\frac{\sqrt{2}}{2},\infty)$, $g(s):=\frac{sf'(s)}{f(s)}=\frac{p-2}{2}\Big[1+\frac{1}{2^{\frac{p-2}{4}}s^{\frac{p-2}{2}}-1}\Big]$ is non-increasing on $(\frac{\sqrt{2}}{2},\infty)$ and $g(s) > \frac{p-2}{2} >-1.$
\end{itemize}

In view of (i)-(ii) the conditions (H1)-(H3) hold and thus \cite[Theorem 2]{SeT} and \cite[Theorem 1.4]{Jeanjean-Zhang-Zhong-2025} applies.

\section{Proofs of Theorem \ref{main-th-1} and Theorem \ref{main-th-2}}

\subsection{Pohozaev type constraint}
For any critical point $u \neq 0$ of $I\big|_{D_a}$, it is standard to show that $u\in \mathcal{P}_a$. Recalling the fiber map $u\mapsto t\star u$,
$(t\star u)(x):=t^{\frac{N}{2}}u(tx),\quad t>0,$
which preserves the $L^2$-norm,  one can check that
$P(u)=\frac{\mathrm{d}}{\mathrm{d}t}I(t\star u)\Big|_{t=1}$.
We decompose $\mathcal{P}=\mathcal{P}_+\cup \mathcal{P}_0\cup \mathcal{P}_-,$
where
$$\mathcal{P}_{+(resp. 0,-)}:=\left\{u\in \mathcal{P}: \frac{\mathrm{d}^2}{\mathrm{d}t^2}I(t\star u)\Big|_{t=1}>0 (resp. =0, <0)\right\},$$
and by a direct calculation
\begin{equation*}
\frac{\mathrm{d}^2}{\mathrm{d}t^2}I(t\star u)\Big|_{t=1}\
=\|\nabla u\|_2^2+(N+2)(N+1)V(u)
-\frac{(p-2)N[(p-2)N-2]}{4p}\|u\|_p^p.
\end{equation*}
For any $u\in \mathcal{P}_a$, we also note that
\beq\label{eq:20240827-we7}
I(u)=\frac{(p-2)N-4}{2(p-2)N}\|\nabla u\|_2^2 +\frac{(p-4)N-4}{(p-2)N}V(u).
\eeq

\bl\lab{lemma:20240826-xl1}
Let $A,B,C>0, \,  \sigma>N+2$, $h(t):=At^2+Bt^{N+2}-C t^{\sigma},\, t>0$. Then $h$ has a unique positive critical point $t_0$ such that $h(t_0)=\max_{t>0}h(t)$
and $h''(t_0)<0$.
\el
\bp
Suppose there exists some $t>0$ such that $h'(t)=0$ and $h''(t)\geq 0$, i.e.,
$$
2At+(N+2)Bt^{N+1}=C\sigma t^{\sigma-1} \mbox{ and }
2A+(N+2)(N+1)Bt^N\geq C\sigma(\sigma-1)t^{\sigma-2}.
$$
By a direct computation,
$A\leq \frac{(N+2-\sigma)\sigma C}{2N}t^{\sigma-2}.$
Noting that $\sigma>N+2$ and $C>0$, we see that $A<0$, a contradiction to the fact $A>0$.

Consequently, $h(t)$ has no local minimizer. Noting that $h(t)>0$ for $t>0$ small, and $h(t)\rightarrow -\infty$ as $t\rightarrow \infty$, we conclude that $\max_{t>0}h(t)$ is achieved at some $t_0>0$. If $h(t)$ has at least two critical points, then it must have a local minimizer, a contradiction.
\ep

\bc\lab{cro:20240829-0911}
Let $N\geq 1$ and $4+\frac{4}{N}<p<2\cdot 2^*$. Then for any $0\neq u\in D_a$, there exists a unique $t=t_u>0$ such that $t\star u\in \mathcal{P}_a$. Furthermore, $t_u<(resp.=,>)1$ if and only if $P(u)<(resp. =,>)0$. In particular, $\mathcal{P}_a = \mathcal{P}_{-} \cap D_a$.
\ec
\bp
We have $t\star u\in D_a$ for all $t>0$. By a direct computation,
$I(t\star u)=\frac{1}{2}\|\nabla u\|_2^2 t^2+V(u) t^{N+2}-\frac{1}{p}\|u\|_p^p t^{\frac{(p-2)N}{2}}.$
Noting that $\frac{(p-2)N}{2}>N+2$ due to $p>4+\frac{4}{N}$, let us define $h(t):=At^2+Bt^{N+2}-Ct^\sigma$ with
$A:=\frac{1}{2}\|\nabla u\|_2^2,$ $ B:=V(u),C:=\frac{1}{p}\|u\|_p^p,$ $\sigma:=\frac{(p-2)N}{2}.$
Then $t\star u\in \mathcal{P}_a$ if and only if $\frac{\mathrm{d}}{\mathrm{d}t}I(t\star u)=0$, i.e., $h'(t)=0$. By Lemma \ref{lemma:20240826-xl1}, we obtain the uniqueness of $t=t_u$. Furthermore, $t_u$ attains the maximum of $h(t)$ in $t\in \R^+$ and
$h'(t)>0$ for $0<t<t_u$ while $h'(t)<0$ for $t>t_u.$
So combining with $P(u)=h'(1)$, we obtain that
$$
P(u)<(resp. =,>)0\Leftrightarrow h'(1)<(resp. =,>)0
\Leftrightarrow t_u<(resp. =,>)1.
$$
\ep

%\br\lab{remark:20240827-1757}
%By Corollary \ref{cro:20240829-0911}, for any $0\neq u\in D_a$, there exists a unique $t=t_u>0$ such that $t\star u \in \mathcal{P}_a$. Hence,
%\beq\lab{eq:20240827-1800}
%M_a:=\inf_{\mathcal{P}_a}I(u)=\inf_{u\in D_a\backslash\{0\}}\max_{t>0}I(t\star u).
%\eeq
%\er

\bl\lab{lemma:20240827-wl2}
Let $N\geq 1$ and $4+\frac{4}{N}<p<2\cdot 2^*$.
There exists $C_{N,p}>0$ such that,
%for any $a>0$,
$$\inf_{u\in \mathcal{P}_a}V(u)\geq C_{N,p}a^{\frac{(N-2)p-4N}{(p-4)N-4}}=:\delta_a>0, \quad \mbox{for any } a >0.$$
\el
\bp
Since $p<2\cdot 2^*=\frac{4N}{(N-2)_+}$, it holds that $\frac{4N-p(N-2)}{2(N+2)}>0$. So for any $u\in D_a$,
recalling Eq.\eqref{eq2.1}, we have
$\|u\|_p^p\leq C(N, p) a^{\frac{4N-p(N-2)}{2(N+2)}} \left(V(u)\right)^{\frac{N(p-2)}{2(N+2)}}.$
Then, for any $u\in \mathcal{P}_a$, it holds that
$$
V(u)< \frac{1}{N+2}\|\nabla u\|_2^2+V(u)
=\frac{(p-2)N}{2p(N+2)}\|u\|_p^p
\leq \frac{(p-2)N}{2p(N+2)} C(N, p) a^{\frac{4N-p(N-2)}{2(N+2)}} \left(V(u)\right)^{\frac{N(p-2)}{2(N+2)}}.
$$
%\begin{align*}
%V(u)< &\frac{1}{N+2}\|\nabla u\|_2^2+V(u)
%=\frac{(p-2)N}{2p(N+2)}\|u\|_p^p\\
%\leq&\frac{(p-2)N}{2p(N+2)} C(N, p) a^{\frac{4N-p(N-2)}{2(N+2)}} \left(V(u)\right)^{\frac{N(p-2)}{2(N+2)}}.
%\end{align*}
Since $p>4+\frac{4}{N}$, one has $\frac{N(p-2)}{2(N+2)}>1$ and it follows from the above inequality that
%\beq\label{eq:20240827-we6}
\begin{align*}
V(u)\geq &\left[\frac{2p(N+2)}{(p-2)N}\frac{1}{C(N, p)} a^{\frac{p(N-2)-4N}{2(N+2)}}\right]^{\frac{2(N+2)}{N(p-2)-2(N+2)}}
=:C_{N,p} a^{\frac{(N-2)p-4N}{(p-4)N-4}}.
\end{align*}
%\eeq
\ep

%\bl\lab{lemma:20240828-1144}
%Let $N\geq 1$ and $4+\frac{4}{N}<p<2\cdot 2^*$.
%Then for any $a>0$, $I\big|_{\mathcal{P}_a}$ is coercive in the sense that
%$\displaystyle \lim_{u\in \mathcal{P}_a, [\|\nabla u\|_2^2+V(u)]\rightarrow +\infty}I(u)=+\infty.$
%\el
%\bp
%It follows by Eq.\eqref{eq:20240827-we7}.
%\ep

For any $u\in \mathcal{P}_a$, let $u^*$ be the Schwarz symmetrization of $u$.
%Then we have the following result.
\bl\lab{lemma:20240829-0941}
Let $N\geq 1$ and $4+\frac{4}{N}<p<2\cdot 2^*$. For any $u\in \mathcal{P}_a$, there exists a unique $t=t_{u^*}\in (0,1]$ such that $t\star u^*\in \mathcal{P}_a$ and $I(t\star u^*)\leq I(u)$.
\el
\bp
%For any $u\in \mathcal{P}_a$, we have that $u\neq 0$ and $P(u)=0$.
By $\|u^*\|_2^2=\|u\|_2^2$, we see that $u^*\in D_a\backslash\{0\}$. Then by Corollary \ref{cro:20240829-0911}, there exists a unique $t=t_{u^*}>0$ such that $t\star u^*\in \mathcal{P}_a$. In particular, by the properties of rearrangement,
$$
\|\nabla u^*\|_2^2\leq \|\nabla u\|_2^2, \,
V(u^*)=\frac{1}{4}\|\nabla (u^*)^2\|_2^2\leq \frac{1}{4}\|\nabla (u^2)\|_2^2=V(u), \,
\|u^*\|_p^p=\|u\|_p^p,
$$
where we have used that $(u^2)^*=(u^*)^2$.
This implies that $P(u^*)\leq P(u)=0$. Hence, by Corollary \ref{cro:20240829-0911} again,  $t_{u^*}\leq 1$. We remark that for any $u\neq 0$, $I(u^*)\leq I(u)$. Combined with the fact that
$s\star u^*=(s\star u)^*~\hbox{for all}~s\in \R^+,$
we finally obtain
\begin{equation}\label{eq:20260410-2001}
I(t\star u^*)=I((t\star u)^*)\leq I(t\star u)
\leq  \max_{s>0}I(s\star u)=I(u),
\end{equation}
proving the lemma.
\ep
\subsection{Properties of $M_a$}
Combining Eq.\eqref{eq:20240827-we7} with Lemma \ref{lemma:20240827-wl2}, one see that
\beq\lab{eq:20240827-1803}
M_a>\frac{(p-4)N-4}{(p-2)N}V(u)\geq \frac{(p-4)N-4}{(p-2)N}C_{N,p}a^{\frac{(N-2)p-4N}{(p-4)N-4}}>0.
\eeq
For any $a_1>a_2>0$, since $D_{a_2}\subsetneqq D_{a_1}$, it follows that $M_a$ is decreasing in $\R^+$. Furthermore, we can establish the following properties for $M_a$.

\bl\lab{lemma:20240830-1313}
Let $N\geq 1$ and $4+\frac{4}{N}<p<2\cdot 2^*$. The function $a \in (0, + \infty) \mapsto M_a$ is continuous.
\el
\bp
Since $M_a$ is decreasing with respect to $a\in \R^+$, it suffices to prove that $M_{a,-}=M_a=M_{a,+}$. Let $a>0$.

Firstly, we consider $a_n\uparrow a$. Clearly $\lim_{n\rightarrow \infty}M_{a_n}\geq M_a$ due to $M_{a_n}\geq M_a, \forall n\in \N$. We only need to prove that $\lim_{n\rightarrow \infty}M_{a_n}\leq M_a+\varepsilon$ for any $\varepsilon>0$.  So, for any $\varepsilon>0$, let $u=u^*\in \mathcal{P}_{a}$ be such that $I(u)\leq M_a+\varepsilon$ and define
$$A:=\frac{1}{2}\|\nabla u\|_2^2,\, B:=V(u), \,  C:=\frac{1}{p}\|u\|_p^p, \, \sigma:=\frac{(p-2)N}{2}.$$
Since $u\in \mathcal{P}_{a}$, we have that
\beq\lab{eq:20240903-0847}
2A+(N+2)B-\sigma C=0.
\eeq
Define $u_n:=\sqrt{\frac{a_n}{a}}u$, then $u_n\in D_{a_n}$.
Let $t_n>0$ be such that $t_n\star u_n\in \mathcal{P}_{a_n}$. Precisely, by Lemma \ref{lemma:20240826-xl1}, $t_n>0$ is uniquely determined by
\beq\lab{eq:20240903-0856}
2A_n+(N+2)B_n t_n^{N}-\sigma C_n t_{n}^{\sigma-2}=0,
\eeq
where
$$\begin{cases}
A_n:=\frac{1}{2}\|\nabla u_n\|_2^2=\frac{1}{2}\frac{a_n}{a}\|\nabla u\|_2^2\rightarrow A,\\
B_n:=V(u_n)=\frac{a_n^2}{a^2}V(u)\rightarrow B,\\
C_N:=\frac{1}{p}\|u_n\|_p^p=(\frac{a_n}{a})^{\frac{p}{2}}C\rightarrow C.
\end{cases}
$$
By Eq.\eqref{eq:20240903-0856} and $\sigma>N+2$, we see that $\{t_n\}$ is bounded. So, passing to a subsequence if necessary, we may assume that $t_n\rightarrow t_0$. Then, by letting $n\rightarrow \infty$ in Eq.\eqref{eq:20240903-0856}, we have
\beq\lab{eq:20240903-0859}
2A+(N+2)B t_{0}^{N}-\sigma C t_{0}^{\sigma-2}=0,
\eeq
which implies that $t_0$ is a critical point of $h(t):=At^2+B t^{N+2}-Ct^\sigma$. So by Lemma \ref{lemma:20240826-xl1} again, combining  Eq.\eqref{eq:20240903-0847} and Eq.\eqref{eq:20240903-0859}, we conclude that $t_0=1$.
Hence,
\begin{align*}
M_{a_n}:=&\inf_{v\in\mathcal{P}_{a_n}}I(v)\leq I\left(t_n\star \sqrt{\frac{a_n}{a}}u\right)
\rightarrow I(u)\leq M_a+\varepsilon,
\end{align*}
proving that $M_{a,-}=M_a$.

Next, we consider $a_n\downarrow a$. Note that $M_{a,+}$ exists and $\lim_{n\rightarrow \infty} M_{a_n}\leq M_a$. Suppose there exists some $\delta>0$ such that
$M_{a_n}<M_a-\delta, \, \forall n\in \N$ and let $u_n=u_n^*\in \mathcal{P}_{a_n}$ be such that
$I(u_n)\leq M_{a_n}+\frac{\delta}{4}.$
Then we have
\begin{equation}\label{L1}
I(u_n)\leq M_a-\frac{\delta}{2}, \forall n\in \N.
\end{equation}
Without loss of generality, since $a_n\downarrow a$,  we may assume that $a_n\leq 2a, \forall n\in \N$.
So, $\{u_n\}\subset \mathcal{P}_{2a}$ with $I(u_n)\leq M_a-\frac{\delta}{2}$ and we deduce from Eq.\eqref{eq:20240827-we7} that $\{u_n\}$ is bounded in $X$. Then by \cite[Proposition 1.7.1]{Cazenave} and \cite[Theorem A.I]{BL}, there exists a subsequence, still denote by $\{u_n\}$, such that $u_n\rightharpoonup u$ in $H^1(\R^N), u_n\rightarrow u$ in $L^q(\R^N)$ with $q\in (2, 2\cdot 2^*)$, and $u_n\rightarrow u$ a.e. in $\R^N$. Noting that $V(u)=\frac{1}{4}\|\nabla (u^2)\|_2^2$, we also have $u_n^2\rightharpoonup u^2$ in $D_{0}^{1,2}(\R^N)$. Then, by the weak lower semi-continuity,
\beq\lab{eq:20240830-1357}
\begin{cases}
\|u\|_2^2 \leq \liminf_{n\rightarrow\infty}\|u_n\|_2^2, & \|\nabla u\|_2^2 \leq \liminf_{n\rightarrow\infty}\|\nabla u_n\|_2^2,\\
V(u)\leq\liminf_{n\rightarrow\infty}V(u_n), & \|u\|_p^p=\liminf_{n\rightarrow\infty}\|u_n\|_p^p.
\end{cases}
\eeq

We claim that $u\neq 0$. If not, by $\|u_n\|_p^p=o_n(1)$ and $P(u_n)=0$, we have that
$\|\nabla u_n\|_2^2+(N+2)V(u_n)=o_n(1),$
which implies $\|\nabla u_n\|_2^2=o_n(1), V(u_n)=o_n(1)$,
and thus contradicts Lemma \ref{lemma:20240827-wl2}.
Then we have from Eq.\eqref{eq:20240830-1357} that $P(u)\leq 0$ and $0<\|u\|_2^2\leq a$ (which implies $u\in D_a\backslash\{0\}$). Recalling Corollary \ref{cro:20240829-0911}, there exists some $0<t\leq1$ such that $t\star u\in \mathcal{P}_a$. So, by Eq.\eqref{eq:20240827-we7} and Eq. \eqref{eq:20240830-1357},
\begin{align*}
M_a=\inf_{v\in \mathcal{P}_a}I(v)\leq I(t\star u) =& \frac{(p-2)N-4}{2(p-2)N}\|\nabla u\|_2^2 t^2 +\frac{(p-4)N-4}{(p-2)N}V(u) t^{N+2}\\
\leq&\frac{(p-2)N-4}{2(p-2)N}\|\nabla u\|_2^2 +\frac{(p-4)N-4}{(p-2)N}V(u)\\
\leq&\liminf_{n\rightarrow \infty}\frac{(p-2)N-4}{2(p-2)N}\|\nabla u_n\|_2^2 +\frac{(p-4)N-4}{(p-2)N}V(u_n)\\
=&\liminf_{n\rightarrow \infty}I(u_n)~\hbox{(since $u_n\in \mathcal{P}_{2a}$)},
\end{align*}
which is a contradiction in view of Eq.\eqref{L1}.
\ep
\br\lab{remark:20240928-1604}
For any $4+\frac{4}{N}<p<2\cdot 2^*$, $M_a$ is continuous, decreasing and $M_a>0$ for any $a>0$. Hence, $\lim_{a\rightarrow \infty}M_a$ exists. We denote this limit by $\tau_p$.
We  see in the following Lemma \ref{lemma:20240902-1541} that $\tau_p=0$ provided $p<2^*$ and $1\leq N\leq 3$.
\er

\bl\lab{lemma:20240902-1541}
Let $1\leq N\leq 3$ and $4+\frac{4}{N}<p<2^*$, then
$M_a\rightarrow 0$ as $a\rightarrow \infty$.
\el
\bp
Let $u=u^*$ with $\|u\|_2^2=1$ and $u\in \mathcal{P}_1$. Define $u_a:=\sqrt{a}u\in \mS_a$. Let
$$A_a:=\frac{1}{2}\|\nabla u_a\|_2^2=\frac{a}{2}\|\nabla u\|_2^2, \quad B_a:=V(u_a)=a^2V(u), \quad C_a:=\frac{1}{p}\|u_a\|_p^p= \frac{a^{\frac{p}{2}}}{p}\|u\|_p^p$$
and
$\displaystyle h_a(t):=A_a t^2+B_at^{N+2}-C_a t^{\sigma}~\hbox{with}~\sigma=\frac{(p-2)N}{2}.$
Recalling Lemma \ref{lemma:20240826-xl1}, there exists an unique critical point $t_a>0$ for $h_a(t)$ in $\R^+$. Precisely,  $t_a>0$ is determined by
$2A_a+(N+2)B_a t_a^{N}-\sigma C_a t_{a}^{\sigma-2}=0$, that
 is,
\beq\lab{eq:20240903-1023}
\begin{aligned}
0=&\|\nabla u\|_2^2+(N+2)V(u) at_a^N- \frac{\sigma}{p}\|u\|_p^p a^{\frac{p-2}{2}} t_{a}^{\sigma-2}\\
=&\|\nabla u\|_2^2+(N+2)V(u) at_a^N- \frac{\sigma}{p}\|u\|_p^p a^{\frac{2}{N}} (at_a^N)^{\frac{\sigma-2}{N}}.
\end{aligned}
\eeq
Since $a\rightarrow \infty$, by Eq.\eqref{eq:20240903-1023}, we see that $at_a^N\rightarrow 0$ and thus
$t_a=O(a^{-\frac{p-2}{2(\sigma-2)}})~\hbox{as}~a\rightarrow \infty.$
So combining  Eq.\eqref{eq:20240827-we7} and Eq. \eqref{eq:20240903-1023}, we have
\begin{align*}
0<M_a=\inf_{v\in \mathcal{P}_a}I(v)\leq I(t_a\star u_a) =&\frac{(p-2)N-4}{2(p-2)N}\|\nabla u\|_2^2 a t_a^2 +\frac{(p-4)N-4}{(p-2)N}V(u) a^2t_{a}^{N+2}\\
\leq&\frac{(p-2)N-4}{2(p-2)N}\frac{\sigma}{p}\|u\|_p^p a^{\frac{p}{2}}t_{a}^{\sigma}+\frac{(p-4)N-4}{(p-2)N}V(u) a^2t_{a}^{N+2}\\
=& O(a^{\frac{p}{2}-\frac{(p-2)\sigma}{2(\sigma-2)}})+O(a^{2-\frac{(N+2)(p-2)}{2(\sigma-2)}})\\
=&O(a^{\frac{p}{2}-\frac{(p-2)\sigma}{2(\sigma-2)}})~ \hbox{due to $\frac{p}{2}-\frac{(p-2)\sigma}{2(\sigma-2)}>2-\frac{(N+2)(p-2)}{2(\sigma-2)}$}.
\end{align*}
Since $\sigma=\frac{(p-2)N}{2}$, we note that
$\frac{p}{2}-\frac{(p-2)\sigma}{2(\sigma-2)}<0$ if and only if $(N-2)p<2N.$

So, $O(a^{\frac{p}{2}-\frac{(p-2)\sigma}{2(\sigma-2)}})\rightarrow 0$ as $a\rightarrow \infty$. Consequently, $M_a\rightarrow 0$ as $a\rightarrow \infty$.
\ep
\bt\lab{th:20240830-1337}
Let $N\geq 1$ and $4+\frac{4}{N}<p<2\cdot 2^*$. Then, for all $a>0$, $M_a$ is achieved by a Schwarz symmetric function.
\et
\bp
By Lemma \ref{lemma:20240829-0941}, we can take a minimizing sequence $\{u_n\}\subset \mathcal{P}_a$ of Schwarz symmetric functions, i.e.,
$u_n=u_n^*,\,  0<\|u_n\|_2^2\leq a,\,  P(u_n)=0, \, I(u_n)\rightarrow M_a,~\hbox{as}~n\rightarrow \infty.$
Eq.\eqref{eq:20240827-we7} implies that $\{\|\nabla u_n\|_2^2\}, \{V(u_n)\}$ are bounded. As in the proof of Lemma \ref{lemma:20240830-1313}, we may assume that $u_n\rightharpoonup u$ in $H^1(\R^N), \, u_n\rightarrow u$ in $L^q(\R^N)$ with $q\in (2, 2\cdot 2^*)$, $u_n\rightarrow u$ a.e. in $\R^N$,  $u_n^2\rightharpoonup u^2$ in $D_{0}^{1,2}(\R^N)$ and Eq.\eqref{eq:20240830-1357} also holds here.
We claim that $u\neq 0$. If not, by $\|u_n\|_p^p=o_n(1)$ and $P(u_n)=0$, we have that
$\|\nabla u_n\|_2^2+(N+2)V(u_n)=o_n(1),$
which implies that $\|\nabla u_n\|_2^2=o_n(1), V(u_n)=o_n(1)$.
So,
$$M_a+o_n(1)=I(u_n)=\frac{1}{2}\|\nabla u_n\|_2^2+V(u_n)-\frac{1}{p}\|u_n\|_p^p=o_n(1),$$
which contradicts the fact $M_a>0$, see Eq.\eqref{eq:20240827-1803}.

Now, by Eq.\eqref{eq:20240830-1357},
$P(u)\leq \liminf_{n\rightarrow \infty} P(u_n)=0.$
Then by Corollary \ref{cro:20240829-0911}, there exists some $t_u\in (0, 1]$ such that $t_u\star u\in \mathcal{P}_a$. At this point reasoning as in
Lemma \ref{lemma:20240830-1313}, we get
\begin{equation*}
M_a=\inf_{v\in \mathcal{P}_a}I(v)\leq I(t_u\star u) \leq \liminf_{n\rightarrow\infty}I(u_n) = M_a,
\end{equation*}
implying that $t_u=1,I(u)=I(t_u\star u)=M_a$, $\|\nabla u_n\|_2^2\rightarrow \|\nabla u\|_2^2$ and $V(u_n)\rightarrow V(u)$  as $n\rightarrow \infty.$
Hence, $u\in \mathcal{P}_a$ attains $M_a$. In particular, $u$ is Schwarz symmetric.
\ep

We now show that any minimizer of $M_a$ has a sign.
%\textcolor{blue}{sign-preserving} critical point of $I$ restricted to $D_a$.

\begin{lemma}\label{lemma:20240830-1318}
Let $N \geq 1$ and $4 + \frac{4}{N} < p < 2 \cdot 2^*$. Suppose $u \in \mathcal{P}_a$ satisfies $I(u) = M_a$. Then there exists $\lambda \in \mathbb{R}$ such that $u$ is a solution to Eq.\,\eqref{eq1.2}.
In particular, up to a sign change, $u$ is positive in $\mathbb{R}^N$.
\end{lemma}

\begin{proof}
Let $b := \|u\|_2^2$, so $0 < b \leq a$. By the monotonicity of $M_a$ in $a$, $M_a = I(u) \geq M_b \geq M_a$, which forces $M_a = M_b$.
Define
$\widetilde{M}_b := \inf_{v \in \mathcal{P} \cap \mathcal{S}_b} I(v).$
Trivially, $\widetilde{M}_b \geq M_b$. Since $u \in \mathcal{S}_b$, we obtain
$M_a = I(u) \geq \widetilde{M}_b \geq M_b = M_a,$
so $I(u) = \widetilde{M}_b$. Applying \cite[Lemma 2.5]{ZCW}, we conclude that $u$ is a solution to Eq.\,\eqref{eq1.2}.

Now $|u| \in \mathcal{P}_a$ also satisfies $I(|u|) = M_a$ and is a solution to Eq.\,\eqref{eq1.2}. Thus $v := \varphi^{-1}(|u|) \in H^1(\mathbb{R}^N)$ is a non-negative solution to Eq.\,\eqref{eq:20240904-0823}. By Lemma \ref{lemma:20260413-2245} and the strict monotonicity of $\varphi$, we conclude that $|u| = \varphi(v) > 0$ in $\mathbb{R}^N$.
\end{proof}
It remains to prove that the  minimizers of $M_a$ belong to  $\mS_a$, and therefore solve the initial problem.
\br\lab{remark:20240830-1333}
%Since, by Lemma \ref{lemma:20240830-1313},  $M_a$ is a continuous decreasing function in $\R^+$: either
Since $M_a$ is a continuous decreasing function in $\R^+$, either (i) or (ii) below must hold,
\begin{itemize}
\item[(i)] there exists $0<a_1<a_2$ such that $M_a\equiv const., \forall a\in [a_1, a_2]$.
\item[(ii)]$M_a$ decreases strictly in $\R^+$.
\end{itemize}
In the next subsection, we prove that only (ii) happens if $1\leq N\leq 4$, see Lemma \ref{th:20240830-1727}.
\er

%\br\lab{remark:20240830-1333}
%Since, by Lemma \ref{lemma:20240830-1313},  $M_a$ is a continuous decreasing function in $\R^+$: either
%Since $M_a$ is a continuous decreasing function in $\R^+$: either
%\begin{itemize}
%\item[(i)] there exists $0<a_1<a_2$ such that $M_a\equiv const., \forall a\in [a_1, a_2]$
%\end{itemize}
%or
%\begin{itemize}
%\item[(ii)]$M_a$ decreases strictly in $\R^+$.
%\end{itemize}
%In the next subsection, we prove that only (ii) happens if $1\leq N\leq 4$, see Lemma \ref{th:20240830-1727}.
%\er

\subsection{Property of the Lagrange multiplier}
\bl\lab{lemma:20240830-1442}
Assume that $N\geq 1, 4+\frac{4}{N}<p<2\cdot 2^*$ and
let $\lambda \in \R$ be given by Lemma \ref{lemma:20240830-1318}. Then $\lambda\geq 0$. In particular, $\lambda=0$ provided $\|u\|_2^2<a$.
\el
\bp
By $0<\|u\|_2^2\leq a$, we can find some $\varepsilon>0$ such that
$$(1+s)u\in D_a,\forall s\in \Lambda, ~\hbox{with}~\Lambda:=\begin{cases}
(-\varepsilon, \varepsilon)~\hbox{if}~\|u\|_2^2<a,\\
(-\varepsilon, 0]~\hbox{if}~\|u\|_2^2=a.
\end{cases}$$
For any $s\in \Lambda$, we consider $I(t\star ((1+s)u))$,
$$I(t\star ((1+s)u))=\frac{1}{2}\|\nabla u\|_2^2 (1+s)^2 t^2 + V(u)(1+s)^4 t^{N+2}-\frac{1}{p}\|u\|_p^p (1+s)^p t^{\frac{(p-2)N}{2}}.$$
Define
$$A(s):=\frac{1}{2}\|\nabla u\|_2^2 (1+s)^2, \,  B(s):=V(u)(1+s)^4, \,  C(s):=\frac{1}{p}\|u\|_p^p (1+s)^p,$$
then by Lemma \ref{lemma:20240826-xl1}, we obtain that: for any $s>0$, there exists a unique $t=t(s)>0$ such that
$$I(t(s)\star ((1+s)u))=\max_{t>0}I(t\star ((1+s)u)).$$
Precisely, $t=t(s)$ is determined by
\beq\lab{eq:20240830-1506}
\|\nabla u\|_2^2(1+s)^2 +(N+2)V(u)(1+s)^4 t^N -\frac{(p-2)N}{2p}\|u\|_p^p(1+s)^p t^{\frac{(p-2)N-4}{2}}=0.
\eeq
In particular, the implicit function theorem implies that $t(s)\in C^2(\R^+)$ and $t(0)=1$.
Define $\rho(s):=I(t(s)\star ((1+s)u)), s\in \Lambda.$
%\textcolor{red}{and the construction of this type of path can also  be found in  \cite[Lemma 8.2]{Jeanjean-Zhang-Zhong-2024} and \cite[Lemma 3.6]{BFG}.}
Since $t(s)\star ((1+s)u)\in \mathcal{P}_a$, we have that $\rho(s)\geq M_a=\rho(0)$. Hence,
\beq\lab{eq:20240830-1713}
\rho'(0)\begin{cases}
=0,\quad&\hbox{if}~\|u\|_2^2<a,~\hbox{since $0$ is an inner point of $\Lambda$},\\
\leq 0,\quad&\hbox{if}~\|u\|_2^2=a,~\hbox{since $0$ is a boundary point of $\Lambda$},
\end{cases}
\eeq
here $\rho'(0)$ means $\rho'_-(0)$ if $\Lambda=(-\varepsilon,0]$.

By direct calculation and using Eq.\eqref{eq:20240830-1506}, we arrive at
\begin{equation*}
%\beq\lab{eq:20240830-1702}
\rho'(s)=\|\nabla u\|_2^2(1+s)t(s)^2+4V(u)(1+s)^3 t(s)^{N+2} -\|u\|_p^p(1+s)^{p-1}t(s)^{\frac{(p-2)N}{2}}.
%\eeq
\end{equation*}
Substituting $s=0$ and $t(0)=1$, it follows that
$\rho'(0)=\|\nabla u\|_2^2+4V(u)-\|u\|_p^p.$
On the other hand, by Lemma \ref{lemma:20240830-1318}, one see that
$
\|\nabla u\|_2^2+4V(u)+\lambda \|u\|_2^2=\|u\|_p^p.$
Thus,
\beq\lab{eq:20240830-1710}
\lambda=-\frac{\rho'(0)}{\|u\|_2^2}.
\eeq
Now, by Eq.\eqref{eq:20240830-1713} and Eq.\eqref{eq:20240830-1710}, we obtain the desired conclusion.
\ep
%As a direct consequence of Lemma \ref{lemma:20240830-1442} we get
\bc\lab{cro:20240830-c1}
Assume that $N\geq 1, 4+\frac{4}{N}<p<2\cdot 2^*$ and
let $\lambda, u$ be given by Lemma \ref{lemma:20240830-1318}. Then $u\in \mS_a$ provided $\lambda\neq 0$.
\ec
\subsection{The case $1\leq N\leq 4$ and proof of Theorem \ref{main-th-1}}
%Firstly we consider the case of $1\leq N\leq 4$.
%With the help of a Liouville type result, we establish the following result.
\bl\lab{th:20240830-1727}
Assume that $1\leq N\leq 4, 4+\frac{4}{N}<p<2\cdot 2^*$ and let $M_a$ be defined in Eq.\eqref{eq:20220930-e1}, then $M_a$ decreases strictly in $\R^+$.
\el
\bp
It is sufficient to exclude the case (i) in Remark \ref{remark:20240830-1333}. We argue by contradiction and suppose there exists some $0<a_1<a_2$ such that $M_{a_1}=M_{a_2}$. By Theorem \ref{th:20240830-1337}, $M_{a_1}$ is attained by some Schwarz symmetric function $u\in D_{a_1}$. Since, $M_{a_1}=M_{a_2}$ and $u\in D_{a_1}\subset D_{a_2}$, $u$ is also a minimizer of $I$ constrained to $\mathcal{P}_{a_2}$. By Lemma \ref{lemma:20240830-1318}, Lemma \ref{lemma:20240830-1442} and $\|u\|_2^2\leq a_1<a_2$, we see that $v = \varphi^{-1}(u) \in H^1(\R^N)$ is a positive solution to Eq.\eqref{eq:20240904-0823} with $\lambda =0$.
Since $1\leq N\leq 4$ and the right hand side of Eq.\eqref{eq:20240904-0823} is non-negative, this is a contradiction to \cite[Lemma A.2]{ikoma2014compactness}.
\ep
{\bf Proof of Theorem \ref{main-th-1}:}  (i) and (ii) follow from Lemma \ref{lemma:20240830-1313} and Lemma \ref{th:20240830-1727} respectively. To prove (iii), first note that, for any $a>0$, by Theorem \ref{th:20240830-1337}, $M_a$ is achieved. Now let $u \in \mathcal{P}_a$ with $I(u) = M_a$ be arbitrary. In particular, by Lemma \ref{th:20240830-1727}, we have that $u\in \mS_a$. By Lemma \ref{lemma:20240830-1318}, up to a sign change, $u>0$ in $\mathbb{R}^N$. Considering $v=\varphi^{-1}(u)$, as in the proof of Lemma \ref{th:20240830-1727}, we deduce that $\lambda>0$. Now, by Lemma \ref{lemma:20260413-2245}, up to a translation, $v=\varphi^{-1}(u)$ is positive, radial and strictly decreasing w.r.t. $r=|x|$. Since $u=\varphi(v)$, it has the same properties.
\hfill $\Box$

%\br\lab{remark:20240830-1810}
%We can also prove Theorem \ref{main-th-1} directly without the strictly decreasing property of $M_a$.  Indeed, one can prove the Liouville's result: Eq.\eqref{eq:20240830-1747} has no non-trivial non-negative solution, see for example \cite[Lemma 2.3]{JLWang}. Then by Theorem \ref{th:20240830-1337} and Lemma \ref{lemma:20240830-1318}, the minimizer $0\neq u\in D_a$ is Schwarz symmetric. Furthermore, $0\leq v=\int_0^u\sqrt{1+2t^2}\,dt \in H^1(\R^N)$ is also a Schwarz symmetric function and it is a non-trivial solution to
%\beq\lab{eq:20240829-1042}
%-\Delta v=-\lambda \varphi(v)\varphi'(v)+\varphi(v)^{p-1}\varphi'(v)\geq 0~\hbox{in}~\R^N.
%\eeq
%Then by \cite[Lemma A.2]{ikoma2014compactness}, we have that $\lambda>0$ and thus it follows by Corollary \ref{cro:20240830-c1} that $u\in \mS_a$.
%\er

\subsection{The case  $N\geq 5$ and proof of Theorem \ref{main-th-2}}
Note that $-\Delta u\geq 0, 0\leq u\in L^2(\R^N)$ is not enough to guarantee $u\equiv 0$ when $N\geq 5$. The situation is different and we obtain the following result.
\bl\lab{th:20241001-1311}
Assume that $N\geq 5, 4+\frac{4}{N}<p<2\cdot 2^*$. Let $M_a$ be defined in Eq.\eqref{eq:20220930-e1} and $a_0$ be given by Eq.\eqref{eq:20241001-1308}, then $M_a$ decreases strictly in $(0,a_0]$.
\el
\bp
Suppose there exist some $0<a_1<a_2\leq a_0$ such that $M_a\equiv M_{a_1}, \forall a\in [a_1,a_2]$. 	As in the proof of Lemma \ref{th:20240830-1727}, we can find some Schwarz symmetric function $0\neq u\in D_{a_1}$ attaining $M_{a_1}$ and solving Eq.\eqref{eq:20240830-1747}. Then by Proposition \ref{remark:20241001-1250}, we have that $u=u_0$. Thus $a_0=\|u_0\|_2^2=\|u\|_2^2\leq a_1<a_2$, a contradiction to $a_2\leq a_0$.
\ep
\begin{lemma}({\bf Existence for $a = a_0$})\label{lemma:20260408-1829}
Let $N\geq 5$, $4+\frac{4}{N}<p<2\cdot 2^*$ and $a_0$ be given by Eq.~\eqref{eq:20241001-1308}. Then $u_0$ attains $M_{a_0}$, i.e., $I(u_0)=M_{a_0}$.
\end{lemma}
\begin{proof}
Since $u_0 \in \mathcal{P}_{a_0}$, $I(u_0) \ge M_{a_0}$. We argue by contradiction and that suppose $I(u_0) > M_{a_0}$. Then, by the monotonicity of $M_a$, it follows that
\begin{equation}\label{eq:20260408-1838}
M_a \le M_{a_0} < I(u_0) \quad \text{for all } a > a_0.
\end{equation}
For any $a > a_0$, Theorem \ref{th:20240830-1337} implies that $M_{a}$ is achieved by a positive Schwarz symmetric function $u_{a} \in D_{a}$. By Lemmas \ref{lemma:20240830-1318} and \ref{lemma:20240830-1442}, $u_a$ solves Eq.\eqref{eq1.2} with some $\lambda_a\geq 0$.
If $\lambda_a = 0$, then by Proposition \ref{remark:20241001-1250}, we would have $u_a = u_0$. Consequently, $M_a = I(u_a) = I(u_0) > M_{a_0}$, which contradicts Eq.\eqref{eq:20260408-1838}. Thus, $\lambda_a > 0$ and $u_a\in \mS_a$ due to Corollary \ref{cro:20240830-c1}.

Now, take a sequence $\{a_n\} \subset \mathbb{R}$ such that $a_0 < a_n \uparrow +\infty$. For each $n \in \mathbb{N}$, the above argument shows that $M_{a_n}$ is achieved by a Schwarz symmetric function $u_n \in \mathcal{S}_{a_n}$ with a corresponding Lagrange multiplier $\lambda_n > 0$.
From the identities
\[
M_{a_n}=I(u_{n})=\frac{(p-2)N-4}{2(p-2)N}\|\nabla u_n\|_2^2 +\frac{(p-4)N-4}{(p-2)N}V(u_n),
\]
and using Eq.\eqref{eq:20260408-1838} we deduce the uniform bounds
\begin{equation*}
%\label{eq:bounds}
\|\nabla u_n\|_2^2 \le \frac{2(p-2)N}{(p-2)N-4}M_{a_0}, \quad
V(u_n)\le \frac{(p-2)N}{(p-4)N-4}M_{a_0}.
\end{equation*}
Since $P(u_n)=0$, it also holds that
\begin{equation}\label{eq:20241002-1406}
\begin{aligned}
\|u_n\|_p^p
&= \frac{2p}{(p-2)N}\Big[\|\nabla u_n\|_2^2+(N+2)V(u_n)\Big] \\
&\le \frac{2p}{(p-2)N}\Bigg[\frac{2(p-2)N}{(p-2)N-4}+(N+2)\frac{(p-2)N}{(p-4)N-4}\Bigg]M_{a_0}.
\end{aligned}
\end{equation}
Used for $u_n$, Eq.\eqref{eq1.2} implies, using the
%The equation for $u_n$,
%\[
%-\Delta u_n - \Delta(u_n^2) u_n + \lambda_n u_n = u_{n}^{p-1} \quad \text{in } \mathbb{R}^N,
%\]
%implies via the
Nehari identity, that
\begin{equation}\label{eq:20241002-1404}
\lambda_n a_n = \lambda_n\|u_n\|_2^2 = \|u_n\|_p^p - \|\nabla u_n\|_2^2 - 4V(u_n) < \|u_n\|_p^p.
\end{equation}
Since $a_n \to \infty$,
%Eqs. \eqref{eq:bounds},
Eq.\eqref{eq:20241002-1406} and Eq.\eqref{eq:20241002-1404} together yield $\lambda_n \to 0^+$.

By \cite[Theorem 1.3-(iii)]{Genoud-Nodari-2024}, the sequence $\{u_n\}$ converges: $u_n \to u_0$ in $D_{0}^{1,2}(\mathbb{R}^N) \cap C^{2}(\mathbb{R}^N)$, where $u_0$ is the unique positive radial solution to Eq.~\eqref{eq:20240830-1747}. Consequently,
\[
\|\nabla u_n\|_2^2 \to \|\nabla u_0\|_2^2, \quad V(u_n) \to V(u_0), \quad \|u_n\|_p^p \to \|u_0\|_p^p.
\]
This convergence implies $I(u_n) \to I(u_0)$. However, we also have $M_{a_0} \ge M_{a_n} = I(u_n)$ by monotonicity. Taking the limit, we obtain the contradiction
\[
M_{a_0} \ge \lim_{n \to \infty} I(u_n) = I(u_0) > M_{a_0}.
\]
Therefore, our assumption $I(u_0) > M_{a_0}$ is false, and we conclude that $I(u_0) = M_{a_0}$.
\end{proof}
\begin{lemma}\label{th:20241001-1404}
Let $N\geq 5$, $4+\frac{4}{N}<p<2\cdot 2^*$ and $a_0$ be given by Eq.\,\eqref{eq:20241001-1308}.
For any $a > 0$, if $M_a \neq M_{a_0}$, then $M_a$ is attained at some $u_a \in \mathcal{S}_a$. Moreover, any energy ground state (up to translations and sign changes) is positive, radial, and strictly decreasing  with respect to $r = |x|$. The associated Lagrange parameter is positive.
\end{lemma}
\begin{proof}
By Theorem\,\ref{th:20240830-1337} and Lemma\,\ref{lemma:20240830-1318}, $M_a$ is attained by a Schwarz symmetric function $u_a \in D_a$, and there exists $\lambda_a \geq 0$ such that $u_a$ satisfies Eq.\,\eqref{eq1.2}. If $u_a \notin \mathcal{S}_a$, then by Lemma\,\ref{lemma:20240830-1442}, $\lambda_a = 0$. Thus $u_a$ solves Eq.\,\eqref{eq:20240830-1747}, and by Proposition\,\ref{remark:20241001-1250}, necessarily $u_a = u_0$. It follows that $M_a = I(u_a) = I(u_0) = M_{a_0}$, contradicting $M_a \neq M_{a_0}$. Therefore, $u_a \in \mathcal{S}_a$. Since $u_a \neq u_0$ and $\lambda_a \geq 0$, Proposition\,\ref{remark:20241001-1250} implies $\lambda_a > 0$.

Now suppose $w_a \in \mathcal{S}_a$ is {\bfseries any} minimizer attaining $M_a$. By Lemma\,\ref{lemma:20240829-0941}, its symmetric-decreasing rearrangement $w_a^* \in \mathcal{S}_a$ also attains $M_a$, with
\begin{equation}\label{eq:20260414-1035}
\|\nabla w_a^*\|_2^2 = \|\nabla w_a\|_2^2, \quad V(w_a^*) = V(w_a), \quad \|w_a^*\|_p^p = \|w_a\|_p^p.
\end{equation}
Applying the above argument to $w_a^*$, we deduce that its associated Lagrange multiplier $\lambda_{w_a}^* $ is positive. From Eq.\,\eqref{eq:20241025-1321} and Eq.\,\eqref{eq:20260414-1035}, the Lagrange multiplier $\lambda_{w_a}$ associated to $w_a$ is equal to $\lambda_{w_a}^*$ and is therefore positive. Then, turning to the dual formulation and using  Lemma\,\ref{lemma:20260413-2245}, it follows that (up to translation and sign change) $w_a$ is positive, radial and strictly radially decreasing.
\end{proof}
In view of Lemmas \ref{lemma:20260408-1829}, \ref{th:20241001-1404} and of the monotonicity of $M_a$, it is well defined that
\beq\lab{eq:20241002-1350}
\tilde{a}^*:=\sup\{a\geq a_0: M_a=M_{a_0}\}.
\eeq
\begin{lemma}\label{lemma:20241002-1349}
Assume that $N\geq 5$ and $4+\frac{4}{N}<p<2\cdot 2^{*}$. Let $\tilde{a}^{*}$ be defined by
Eq.\eqref{eq:20241002-1350}. Then $\tilde{a}^{*}=+\infty$; that is, $M_{a}\equiv M_{a_{0}}$ for all
$a\in[a_{0},+\infty)$.
\end{lemma}
\begin{proof}
We argue by contradiction and suppose that $\tilde{a}^{*}<\infty$. By Def.\eqref{eq:20241002-1350} and the monotonicity of $a\mapsto M_{a}$, we have
\begin{equation}\label{eq:20241002-1352}
M_{a}<M_{a_{0}}\qquad\text{for all }a>\tilde{a}^{*}.
\end{equation}
Choose a sequence $\{a_{n}\}\subset\mathbb{R}$ with $\tilde{a}^{*}<a_{n}\uparrow +\infty$.
For each $n\in\mathbb{N}$, Lemma \ref{th:20241001-1404} guarantees the existence of a
Schwarz symmetric minimizer $u_{n}\in\mathcal{S}_{a_{n}}$ of $M_{a_{n}}$ and a Lagrange
multiplier $\lambda_{n}>0$ satisfying
\[
-\Delta u_{n}-\Delta(u_{n}^{2})u_{n}+\lambda_{n}u_{n}=u_{n}^{\,p-1}\quad\text{in }\mathbb{R}^{N}.
\]
Following the same line of reasoning as in the proof of Lemma
\ref{lemma:20260408-1829} (specifically, the part establishing $\lambda_{n}\to0^{+}$ and
the consequent convergence $u_{n}\to u_{0}$), we obtain
\[
M_{a_{n}}=I(u_{n})\rightarrow I(u_{0})=M_{a_{0}}.
\]
Because $M_{a}$ is decreasing and $a_{n}>\tilde{a}^{*}$ for all $n$,
%the limit $\lim M_{a_{n}}=M_{a_{0}}$
$M_{a_{n}} \to M_{a_{0}}$ contradicts Eq.\eqref{eq:20241002-1352}.
 %the strict inequality$M_{a_{n}}<M_{a_{0}}$ required by \eqref{eq:20241002-1352}.
Hence our assumption
$\tilde{a}^{*}<\infty$ is false, and we conclude that $\tilde{a}^{*}=+\infty$.
\end{proof}

\bl\lab{th:20241002-1424}
Let $N\geq 5$ and $4+\frac{4}{N}<p<2\cdot 2^*$. Then there is no $u\in \mS_a$ which achieve  $M_a$ provided $a\in (a_0, +\infty)$.
\el
\bp
Suppose that, for a $a>a_0$, some $u_a\in \mS_a$ achieves $M_a$. Then $u_a^*\in \mS_a$ also achieves $M_a$.
 %\red{In view of Lemma \ref{lemma:20240829-0941}, it is not restrictive to assume that $u_a$ is a Schwartz symmetric function.}
%$u_a^* \in \mS_a$ of $u_a$ also achieves $M_a$.}
Then by Lemma \ref{lemma:20240830-1318}, there exists some $\lambda_a\in \R$ such that
$u_a^*$ solves Eq.\eqref{eq1.2}.

However, letting $b>a$, then by Lemma \ref{lemma:20241002-1349}, $M_b=M_a=M_{a_0}$.
Thus, $u_a^* \in \mS_a\subset D_b$ also achieves $M_b$. So, by Lemma \ref{lemma:20240830-1442}, $\lambda_a=0$. Consequently, in view of Proposition \ref{remark:20241001-1250}, $u_a^*=u_0$ and therefore $a=\|u_a^*\|_2^2=\|u_0\|_2^2=a_0$, which is a contradiction.
\ep

\bl\lab{th:20241004-1344}({\bf Uniqueness of energy ground state for $a = a_0$})
Let $N\geq 5,4+\frac{4}{N}<p<2\cdot 2^*$ and $u_0$ be the unique positive radial solution to Eq.\eqref{eq:20240830-1747}.
Then $u_0$ is, up to translation and sign change, the unique energy ground state when $a=a_0$.
\el
\bp
Assume that $u_1\in \mS_{a_0}$ achieves $M_{a_0}$. By Lemma \ref{lemma:20240830-1318}, there exists some $\lambda_1\in \R$ such that $u_1$ solves Eq.\eqref{eq1.2} with $\lambda=\lambda_1$.
For any fixed $a_1>a_0$, by Lemma \ref{lemma:20241002-1349}, $M_{a_1} = M_{a_0}$ and then $u_1$ is also a minimizer of Problem \eqref{eq:20220930-e1} with $a=a_1$. Noting that $\|u_1\|_2^2=a_0<a_1$, by Lemma \ref{lemma:20240830-1442}, $\lambda_1=0$. Hence, if $u_1$ is positive and radial necessarily $u_1=u_0$. In the Appendix we establish that any energy ground state for $a=a_0$ is positive and radial, up to translation and sign change.
\ep
%\red{
%\begin{remark}\label{remark:20260414-1128}
%One may wonder if (up to translation and sign change) the energy ground state is unique when $a=a_0$. Note that if $u_1 \in S_{a_0}$ is a minimizer of $M_{a_0}$ then $u_1^*$ is also an energy ground state and thus $u_1^* = u_0$ by Lemma \ref{th:20241004-1344}. If we could conclude that, up to translation and sign change, $u_1 = u_1^*$ the uniqueness would hold. In that direction note that the associated Lagrange multiplier satisfies $\lambda_{u_1} = 0$ (by proof of Lemma \ref{th:20241004-1344}), that we may assume $u_1 > 0$ (by Lemma \ref{lemma:20240830-1318}) and also that $||\nabla u_1||_2 =||\nabla u_1^*||_2$. However, this information alone is not sufficient to conclude that $u_1 = u_1^*$. The uniqueness issue remains open at this stage.
%\end{remark}
%}
{\bf Proof of Theorem \ref{main-th-2}:}
 (i) follows from Lemma \ref{lemma:20240830-1313} and (ii) follows from Lemma \ref{th:20241001-1311} and Lemma \ref{lemma:20241002-1349}. By Lemma \ref{th:20241001-1311},  $M_a > M_{a_0}$ for any $a < a_0$, thus (iii) follows from Lemma \ref{th:20241001-1404}. Now (iv) follows from Lemma \ref{lemma:20260408-1829} and Lemma \ref{th:20241004-1344}. Finally, (v) is Lemma \ref{th:20241002-1424}. \hfill $\Box$

\s{Asymptotic behavior of the energy ground states}
\subsection{Some preparations}
For any $a \in (0, a_0)$,  let $u_a$ be an energy ground state given by Theorem \ref{main-th-1} or Theorem \ref{main-th-2}. In particular
$u_a=u_a^*\in \mathcal{P}_a$ and  it satisfies
\beq\lab{eq:20240903-1648}
-\Delta u_a-\Delta(u_a^2) u_a +\lambda_a u_a=u_{a}^{p-1}~\hbox{in}~\R^N,
\eeq
with $I(u_n)=M_a$
and  $\lambda_a >0$.
\bl\lab{lemma:20240903-1650}
Let $\lambda_a, u_a$ be given above, then the following properties hold:
\begin{itemize}
\item[(i)] For $N\geq 1$ and $4+\frac{4}{N}<p<2\cdot 2^*$,
$$\|\nabla u_a\|_2^2+V(u_a)\rightarrow \infty \quad \mbox{and} \quad \|u_a\|_p\rightarrow \infty, \quad \mbox{as} \quad  a\rightarrow 0^+.$$

\item[(ii)] For $1\leq N\leq 4$ and $p\in (4+\frac{4}{N},2\cdot 2^*)$, let $\tau_p$ be defined in Remark \ref{remark:20240928-1604}, then
     \beq\lab{eq:20240928-1631}
\tau_{1,p}:=\limsup_{a\rightarrow +\infty}\|\nabla u_a\|_2^2 \leq \frac{2(p-2)N}{(p-2)N-4}\tau_p,
\eeq
\beq\lab{eq:20240928-1633}
\tau_{2,p}:=\limsup_{a\rightarrow +\infty}V(u_a)\leq \frac{(p-2)N}{(p-4)N-4}\tau_p
\eeq
 and
 \beq\lab{eq:20240928-1635}
\tau_{3,p}:=\limsup_{a\rightarrow +\infty}\|u_a\|_p^p\leq \frac{2p}{(p-2)N}\Big[\tau_{1,p}+(N+2)\tau_{2,p}\Big].
\eeq
In particular, $\tau_{1,p}=\tau_{2,p}=\tau_{3,p}=0$ when $1\leq N\leq 3$ and $4+\frac{4}{N}<p<2^*$.
\end{itemize}
\el
\bp
First, recalling Eq.\eqref{eq:20240827-we7}, we observe that
\beq\lab{eq:20240903-1700}
M_a=I(u_a)=\frac{(p-2)N-4}{2(p-2)N}\|\nabla u_a\|_2^2 +\frac{(p-4)N-4}{(p-2)N}V(u_a).
\eeq

(i) Since  $\frac{(N-2)p-4N}{(p-4)N-4}<0$, from Eq.\eqref{eq:20240827-1803} and Eq.\eqref{eq:20240903-1700}, we conclude that
$\displaystyle \|\nabla u_a\|_2^2+V(u_a)\rightarrow \infty, \, \mbox{as} \, a\rightarrow 0^+.$
Now since $P(u_a)=0$, it holds that
\beq\lab{eq:20241002-2052}
\|u_a\|_p^p=\frac{2p}{(p-2)N}\big[\|\nabla u_a\|_2^2+(N+2)V(u_a)\big].
\eeq
Thus also $\|u_a\|_p\rightarrow \infty$ as $a\rightarrow 0^+$.

(ii) As in the proof of Lemma \ref{lemma:20241002-1349}, by Eq.\eqref{eq:20240903-1700}  we observe that
\begin{equation*}
%\beq\lab{eq:20240928-1628}
\|\nabla u_a\|_2^2 \leq \frac{2(p-2)N}{(p-2)N-4}M_a, \quad V(u_a)\leq \frac{(p-2)N}{(p-4)N-4}M_a.
%\eeq
\end{equation*}
This proves Eq.\eqref{eq:20240928-1631} and Eq.\eqref{eq:20240928-1633}.
Recalling Eq.\eqref{eq:20241002-2052}, we conclude that Eq.\eqref{eq:20240928-1635} is also true.
In particular, when $1\leq N\leq 3$ and $4+\frac{4}{N}<p<2^*$, we have $\tau_p=0$ due to Lemma \ref{lemma:20240902-1541}. Thus $\tau_{1,p}=\tau_{2,p}=\tau_{3,p}=0$.
\ep

In Corollary \ref{remark:20240928-2218} we shall determine the values of $\tau_{1,p}, \tau_{2,p}$ and $\tau_{3,p}$ in the remaining cases. Namely, when $p \in [6,12)$ if $N=3$ and when $p \in (5,8)$ if N=4.

%\textcolor{red}{
%\br\lab{remark:20240928-1646}
%In
%In the following, we devote to study the asymptotic behavior of $u_a$ as $a\rightarrow +\infty$. We shall make clear the values of $\tau_{1,p}, \tau_{2,p}$ and $\tau_{3,p}$ for the remaining cases when $1\leq N\leq 4$, i.e., $6\leq p<12, N=3$ and $5<N<8, N=4$, see Remark \ref{remark:20240928-2218}.
%\er
%}

\bl\lab{cro:20240904-0800}
Let $4+\frac{4}{N}<p<2\cdot 2^*$.
\begin{itemize}
\item[(i)] If $1\leq N\leq 4$, then $\lambda_a\rightarrow 0$ as $a\rightarrow +\infty$.
\item[(ii)] If $N\geq 5$, let $a_0$ be given by Eq.\eqref{eq:20241001-1308}, then $\lambda_a\rightarrow 0$ as $a\uparrow a_0$.
\end{itemize}
In particular,
$\lambda_a a\rightarrow 0$ as $a\rightarrow +\infty$ when $1\leq N\leq 3$.
\el
\bp
By Eq.\eqref{eq:20240903-1648}, we first note that
\beq\lab{eq:20240903-1704}
\lambda_a a=\lambda_a\|u_a\|_2^2=\|u_a\|_p^p-\|\nabla u_a\|_2^2-4V(u_a)<\|u_a\|_p^p.
\eeq
(i) If $1\leq N\leq 4$, then, by Eq.\eqref{eq:20240928-1635}
\beq\lab{eq:20240928-1653}
\limsup_{a\rightarrow +\infty}\lambda_a a\leq \tau_{3,p}.
\eeq
Since $\lambda_a>0$ and $a\rightarrow +\infty$, we conclude that $\lambda_a\rightarrow 0$.
Now, when $1\leq N\leq 3$, we know by Lemma \ref{lemma:20240903-1650} that $\tau_{3,p}=0$. Thus from Eq.\eqref{eq:20240928-1653}, we obtain that $\lambda_a a\rightarrow 0$ as $a\rightarrow \infty$.

(ii) Let $N\geq 5$ and consider $a>\frac{1}{2}a_0$.
%By the monotonicity of $M_a$, $M_a\leq M_{\frac{a_0}{2}}$.
Noting that
\begin{equation*}
%\beq\lab{eq:20241002-1847}
\begin{aligned}
\|u_a\|_p^p=&\frac{2p}{(p-2)N}\Big[\|\nabla u_a\|_2^2+(N+2)V(u_a)\Big]\\
<&\frac{2p}{(p-2)N}\Big[\frac{2(p-2)N}{(p-2)N-4}+(N+2)\frac{(p-2)N}{(p-4)N-4}\Big]M_{a},
\end{aligned}
%\eeq
\end{equation*}
by Eq.\eqref{eq:20240903-1704} and $M_a\leq M_{\frac{a_0}{2}}$, one see that $\{\lambda_a a\}_{a\geq \frac{a_0}{2}}$ is bounded.
Consequently, for any $a_n\uparrow a_0$, $\{\lambda_{a_n}\}$ is bounded. Passing to a subsequence, we denote
\begin{equation*}
%\beq\lab{eq:20241002-1832}
\lim_{n\rightarrow \infty}\lambda_{a_n}=:\lambda^*\in [0, +\infty).
%\eeq
\end{equation*}
If $\lambda^*\neq 0$, then $\lambda^*>0$. Without loss of generality, we may assume that $\lambda_{a_n}\ge\frac{1}{2}\lambda^*$ for any $n$. Noting that $(\lambda_{a_n}, u_{a_n})$ satisfies Eq.\eqref{eq:20240903-1648} and $\{\|\nabla u_{a_n}\|_2\}, \{V(u_{a_n})\}, \{\|u_{a_n}\|_2\}$ are bounded, one get
$$
-\Delta v_n+\lambda_{a_n} \varphi(v_n)\varphi'(v_n)=\varphi(v_n)^{p-1}\varphi'(v_n)~\hbox{in}~\R^N,
$$
and that $\{v_n\}$ is bounded in $H_{rad}^1(\R^N)$, where $v_n=\varphi^{-1}(u_{a_n})$. Assume that $v_n\rightarrow v_\infty$ weakly in $H_{rad}^1(\R^N)$ and a.e. in $\R^N$ as $n\rightarrow\infty$, then $v_\infty\in H_{rad}^1(\R^N)\cap C^2(\R^N)$ and satisfies
$$
-\Delta v_\infty+\lambda^\ast \varphi(v_\infty)\varphi'(v_\infty)=\varphi(v_\infty)^{p-1}\varphi'(v_\infty)~\hbox{in}~\R^N.
$$
It follows from \cite[Radial Lemma 1]{Strauss} that $\lim_{|x|\rightarrow\infty}v_n(x)=0$ uniformly for $n$. In particular, since $|\varphi(s)| \leq |s|$, see Eq.\eqref{eq3.9}, taking $R >0$ large enough, for all $n \in \N$,
$$|v_n(x)| \leq 1 \quad \mbox{and} \quad |\varphi(v_n)|^{p-2} \leq \frac{1}{4}\lambda^*, \,\, \forall |x| > R.$$
Also, since  $\varphi'(s)>0$ for $s\geq 0$ and $\varphi'(0)=1$, see Eq.\eqref{eq3.8}, we can find some $C>0$ such that $\varphi(s)\varphi'(s)\geq Cs$ for $s\in [0,1]$.
Thus,
\begin{equation*}
\begin{aligned}
0=&-\Delta v_n+\left[\lambda_n-\varphi(v_n)^{p-2}\right]\varphi(v_n)\varphi'(v_n) \geq -\Delta v_n+ \frac{1}{4} \lambda^*\varphi(v_n)\varphi'(v_n) \\
\geq & -\Delta v_n+ \frac{1}{4} \lambda^* C v_n, \quad \forall |x| > R.
\end{aligned}
\end{equation*}
So, applying a comparison principle, we deduce that $v_n$ decays exponentially at infinity uniformly for $n$.

On the other hand, the Moser iteration technique shows that $\{v_n\}$ is bounded in $L^\infty(\R^N)$, and then standard elliptic regularity implies that $v_n\rightarrow v_\infty$ in $C_{\text{loc}}^{2}(\R^N)$.

 Finally, it follows that $v_n\rightarrow v_\infty$ strongly in $H_{rad}^1(\R^N)\cap C_{\text{loc}}^{2}(\R^N)$ as $n\rightarrow\infty$. Let $u_\infty=\varphi(v_\infty)$, then
\begin{equation*}
%\beq\lab{eq:20241004-1402}
-\Delta u_\infty-\Delta(u_\infty^2) u_\infty +\lambda^\ast u_\infty=u_{\infty}^{p-1}~\hbox{in}~\R^N.
%\eeq
\end{equation*}
The exponential decay of $v_n$ implies the one of $u_{a_n}$. In particular, $\|u_{a_n}\|_2^2\rightarrow \|u_\infty\|_2^2$ and $V(u_{a_n})\rightarrow V(u_\infty)$ as $n\rightarrow \infty$. Thus it follows that $I(u_{a_n})\rightarrow I(u_\infty)$ as $n\rightarrow \infty$. Hence, $u_\infty\in \mS_{a_0}$ with
$\displaystyle
I(u_\infty)=\lim_{n\rightarrow\infty}I(u_{a_n})=\lim_{n\rightarrow\infty}M_{a_n}=M_{a_0}.
$
Then, by
%the uniqueness in
Lemma \ref{th:20241004-1344}, $u_\infty=u_0$, a contradiction to $\lambda^*\neq 0$.
\ep

\bl\lab{lemma:20240904-0808}
Let $N\geq 1$ and $4+\frac{4}{N}<p<2\cdot 2^*$. Then $\|u_a\|_\infty=u_a(0)\rightarrow \infty$ as $a\rightarrow 0^+$.
\el
\bp
We proceed by contradiction, assuming there exists a sequence $a_n\rightarrow 0^+$ and $M>0$ such that
$\sup_{n\in \N}u_{a_n}(0)\leq M.$ Then $\|u_{a_n}\|_p^p\leq \|u_{a_n}\|_{\infty}^{p-2} \|u_{a_n}\|_2^2 \leq M^{p-2}a_n\rightarrow 0,$
which is a contradiction to Lemma \ref{lemma:20240903-1650}-(i).
\ep
In order to get further estimates on the $L^\infty$-norm, we use again the dual approach.
\bl\lab{lemma:20240904-0824}
Let $N\geq 1$.
For any $M>0$, there exists $C_M>0$ such that, for any positive radial decreasing function $v\in H^1(\R^N)$, which solves Eq.\eqref{eq:20240904-0823} with $\lambda\in (0,M)$,
$$\max_{x\in \R^N}v(x)=\|v\|_\infty=v(0)\leq C_M, \forall ~\lambda\in (0,M).$$
\el
\bp
Proceeding by contradiction, we assume there exists a sequence $\{(\lambda_n, v_n)\}\subset (0, M)\times H^1(\R^N)$ where $v_n$ is a positive radial decreasing solution to Eq.\eqref{eq:20240904-0823} with $\lambda=\lambda_n$ and
$\displaystyle M_n:=v_n(0)\rightarrow +\infty~\hbox{as}~n\rightarrow \infty.$

Following a blow up procedure introduced by Gidas and Spruck \cite{GidasSpruck1981}, we perform a rescaling  and define $\tilde{v}_n:\R^N\rightarrow \R$ by
$\displaystyle \tilde{v}_n(y):=\frac{1}{M_n}v_n\left(\frac{y}{M_{n}^{\frac{p-4}{4}}}\right).$
Then $\tilde{v}_n(0)=\max_{y\in \R^N}\tilde{v}_n(y)= \|\tilde{v}_n\|_\infty =1$ and
\beq\lab{eq:20240904-1027}
-\Delta \tilde{v}_n=\frac{1}{M_{n}^{\frac{p-2}{2}}} \left[\varphi(M_n\tilde{v}_n)^{p-1}\varphi'(M_n\tilde{v}_n)-\lambda_n \varphi(M_n\tilde{v}_n)\varphi'(M_n\tilde{v}_n)\right].
\eeq
By Eq.\eqref{eq3.9}, we have
\begin{align*}
&\frac{1}{M_{n}^{\frac{p-2}{2}}} \varphi(M_n\tilde{v}_n)^{p-1}\varphi'(M_n\tilde{v}_n)
=\frac{\varphi(M_n\tilde{v}_n)^{p-1}\varphi'(M_n\tilde{v}_n)M_n\tilde{v}_n }{M_{n}^{\frac{p-2}{2}}M_n\tilde{v}_n}\\
\leq& \frac{\varphi(M_n\tilde{v}_n)^{p-2}\varphi^2(M_n\tilde{v}_n) }{M_{n}^{\frac{p-2}{2}}M_n\tilde{v}_n}
=\frac{\varphi^p(M_n\tilde{v}_n)}{M_{n}^{\frac{p}{2}}\tilde{v}_n}
\leq\frac{[2^{\frac{1}{4}}\sqrt{M_n\tilde{v}_n}]^p}{M_{n}^{\frac{p}{2}}\tilde{v}_n} = 2^{\frac{p}{4}}\tilde{v}_n^{\frac{p-2}{2}}
\leq 2^{\frac{p}{4}}
\end{align*}
and, since $\varphi(s)\varphi'(s)\leq 2^{-\frac{1}{2}}, \forall s\in \mathbb{R}^+$,
\begin{align*}
\frac{1}{M_{n}^{\frac{p-2}{2}}}\lambda_n \varphi(M_n\tilde{v}_n)\varphi'(M_n\tilde{v}_n)
\leq 2^{-\frac{1}{2}}\frac{1}{M_{n}^{\frac{p-2}{2}}}\lambda_n
\rightarrow 0.
\end{align*}
Hence, the RHS of Eq.\eqref{eq:20240904-1027} is in $L^\infty(\R^N)$. So, applying standard elliptic estimates, and passing to a subsequence if necessary, we may assume that $\tilde{v}_n\rightarrow \tilde{v}$ in $C_{\text{loc}}^{1,\alpha}(\R^N)$ where $\tilde{v}$ is radial, bounded, decreasing and non-negative.

Denote $\Omega:=\text{supp}~\tilde{v}$,
then $\Omega=\R^N$ or $\Omega=B_R(0)$ for some $R>0$. In particular, the Schauder estimation implies that $\tilde{v}_n\rightarrow \tilde{v}$ in $C_{loc}^{2}(\Omega)$.

For any $\Omega_1\subset\subset \Omega$, there exists some $C_{\Omega_1}>0$ such that $\inf_{x\in \Omega_1}\tilde{v}_n\geq C_{\Omega_1}, \forall n\in \N$. Then by $M_n\rightarrow +\infty$, one see that $M_n\tilde{v}_n\rightarrow +\infty$ as $n\rightarrow +\infty$ uniformly in $\Omega_1$.

Thus, if $\Omega=\R^N$, then by Eq.\eqref{eq3.8} and Eq.\eqref{eq:20240904-1055}, $\tilde{v}$ is a non-trivial positive bounded radial solution to
$$-\Delta\tilde{v}=2^{\frac{p-4}{4}}\tilde{v}^{\frac{p-2}{2}}~\hbox{in}~\R^N.$$
Since $\frac{p}{2}-1<2^*-1$, this is a contradiction to \cite[Theorem 2.2]{Jeanjean2024}.

If $\Omega=B_R(0)$, since $\tilde{v}\in C^{1,\alpha}_{\text{loc}}(\R^N)$ and $\min_{x\in \R^N}\tilde{v}(x)=0=\tilde{v}$ on $\partial B_R(0)$, we conclude that $\frac{\partial \tilde{v}}{\partial \nu}=0$ on $\partial B_R(0)$. Hence, $\tilde{v}$ is a solution to
$$\begin{cases}
-\Delta\tilde{v}= 2^{\frac{p-4}{4}}\tilde{v}^{\frac{p}{2}-1}~\hbox{in}~B_R(0),\\
\tilde{v}>0~\hbox{in}~B_R(0), \tilde{v}=\frac{\partial \tilde{v}}{\partial \nu}=0~\hbox{on}~\partial B_R(0).
\end{cases}$$
However, by the Hopf's Lemma, see \cite{Evans}, we have that $\frac{\partial \tilde{v}}{\partial \nu}<0$ on $\partial B_R(0)$, which is a contradiction.
\ep

%Now, for any $a>0$, we denote $u_a=\varphi(v_a)$ in the following.

\bc\lab{cro:20240904-1324}
Let $N\geq 1$ and $4+\frac{4}{N}<p<2\cdot 2^*$. Then
$\lambda_a\rightarrow +\infty$ as  $a\rightarrow 0^+$.
\ec
\bp
By Lemma \ref{lemma:20240904-0808},  $\|u_a\|_\infty\rightarrow \infty$ as $a\rightarrow 0^+$ and since $u_a=\varphi(v_a)$, also $\|v_a\|_\infty\rightarrow \infty$ as $a\rightarrow 0^+$. Hence, by Lemma \ref{lemma:20240904-0824}, we conclude that $\lambda_a\rightarrow +\infty$ as  $a\rightarrow 0^+$.
\ep

\subsection{The case $a\rightarrow 0^+$ and proof of Theorem \ref{th:20240905-2122}}
Let $N\geq 1$ and $4+\frac{4}{N}<p<2\cdot 2^*$. By Lemma \ref{lemma:20240904-0808} and Corollary \ref{cro:20240904-1324}, we know that, as $a\rightarrow 0^+$,
$$\lambda_a\rightarrow +\infty, \, \|u_a\|_\infty\rightarrow \infty, \, \|v_a\|_\infty\rightarrow \infty.$$
For $a_n\rightarrow 0^+$, we denote $\lambda_{a_n}, u_{a_n}, v_{a_n}$ by $\lambda_n, u_n, v_n$ respectively for simplicity.

\bl\lab{lemma:20240904-1639}
Let $N\geq 1$ and $4+\frac{4}{N}<p<2\cdot 2^*$. Let $\lambda_n, u_n, v_n$ be defined above. Then
\beq\lab{eq:20240907-1841}
0<\liminf_{n\rightarrow \infty}\frac{\|v_n\|_{\infty}^{\frac{p-2}{2}}}{\lambda_n}\leq \limsup_{n\rightarrow \infty}\frac{\|v_n\|_{\infty}^{\frac{p-2}{2}}}{\lambda_n}<+\infty.
\eeq
%$$\liminf_{n\rightarrow \infty}\frac{\|v_n\|_{\infty}^{\frac{p-2}{2}}}{\lambda_n}>0.$$
\el
\bp
From Eq.\eqref{eq:20240903-1648} we deduce that $\lambda_n \leq ||u_n||_{\infty}^{p-2}$. Observing, see Eq.\eqref{eq3.8}, that
$\|u_n\|_{\infty}=u_n(0)=\varphi(v_n(0))=\varphi(\|v_n\|_{\infty})=(2^{\frac{1}{4}}+o_n(1))\|v_n\|_{\infty}^{\frac{1}{2}},$
we conclude that the LHS of Eq.\eqref{eq:20240907-1841} holds true. For the RHS we argue by contradiction and assume there exists a sequence $a_n\rightarrow 0^+$ such that
\begin{equation*}
%\beq\lab{eq:20240904-1623}
\frac{\|v_n\|_{\infty}^{\frac{p-2}{2}}}{\lambda_n}\rightarrow +\infty.
%\eeq
\end{equation*}
Let $M_n:=\|v_n\|_\infty$ and set $\tilde{v}_n(x):=\frac{1}{M_n}v_n\left(\frac{x}{M_{n}^{\frac{p-4}{4}}}\right)$.
Then $\tilde{v}_n(0)=\|\tilde{v}_n\|_\infty=1$,  Eq. \eqref{eq:20240904-1027} holds and we conclude as in the proof of Lemma \ref{lemma:20240904-0824}.
\ep

{\bf Proof of Theorem \ref{th:20240905-2122}:}
Define
$\displaystyle
\tilde{v}_n(x):=\lambda_{n}^{-\frac{2}{p-2}}v_n\left(\frac{x}{\lambda_{n}^{\frac{p-4}{2(p-2)}}}\right),
$
then by Eq. \eqref{eq:20240904-0823},
%Eq.\eqref{eq:20240907-1834}
a direct computation shows that
\beq\lab{eq:20240907-1840}
-\Delta \tilde{v}_n=-\varphi(\lambda_{n}^{\frac{2}{p-2}}\tilde{v}_n)\varphi'(\lambda_{n}^{\frac{2}{p-2}}\tilde{v}_n)
+\lambda_{n}^{-1}\varphi(\lambda_{n}^{\frac{2}{p-2}}\tilde{v}_n)^{p-1}\varphi'(\lambda_{n}^{\frac{2}{p-2}}\tilde{v}_n)~\hbox{in}~\R^N.
\eeq
Recalling Eq.\eqref{eq:20240907-1841}, $\{\tilde{v}_n\}$ is bounded in $L^\infty(\R^N)$.
Since $\varphi(s)\varphi'(s)\leq 2^{-\frac{1}{2}}, \forall s\in \mathbb{R}^+$,
we have that
$$\varphi(\lambda_{n}^{\frac{2}{p-2}}\tilde{v}_n)\varphi'(\lambda_{n}^{\frac{2}{p-2}}\tilde{v}_n)\leq \frac{\sqrt{2}}{2}, \, \forall x\in \R^N, \, \forall n\in \N.$$
By Eq.\eqref{eq3.9}, we deduce
$$\lambda_{n}^{-1}\varphi(\lambda_{n}^{\frac{2}{p-2}}\tilde{v}_n)^{p-1}\varphi'(\lambda_{n}^{\frac{2}{p-2}}\tilde{v}_n)
\leq  \frac{\sqrt{2}}{2}\lambda_{n}^{-1} \left[2^{\frac{1}{4}}\sqrt{\lambda_{n}^{\frac{2}{p-2}}\tilde{v}_n}\right]^{p-2}
= 2^{\frac{p-4}{4}} \tilde{v}_{n}^{\frac{p-2}{2}} \in L^\infty(\R^N).$$
\noindent
Hence, the RHS of Eq.\eqref{eq:20240907-1840} is in $L^\infty(\R^N)$. Now, as in the proof of Lemma \ref{lemma:20240904-0824},
$$\tilde{v}_n\rightarrow \tilde{v}\in C_{\text{loc}}^{1,\alpha}(\R^N)\cap C^{2}_{\text{loc}}(\Omega) \, \mbox{ with } \,  \Omega:=\text{supp}~\tilde{v},$$
where $\Omega=\R^N$ or  $\Omega=B_R(0)$  for some  $R>0.$ Note that the lower bound in Eq.\eqref{eq:20240907-1841} guarantees that $\tilde{v} \not \equiv 0$.

% $\tilde{v}_n\rightarrow \tilde{v}\in C_{\text{loc}}^{1,\alpha}(\R^N)\cap C^{2}_{\text{loc}}(\Omega)$ with $\Omega:=\text{supp}~\tilde{v}$, and $\Omega=\R^N$ or $\Omega=B_R(0)$ for some $R>0$.

\noindent
For $x\in \Omega$, since $\tilde{v_n}(x) \to \tilde{v}(x) >0$, we have, using Eq.\eqref{eq3.8}  and  Eq.\eqref{eq:20240904-1055}, that
\begin{equation*}
%\beq\lab{eq:20240907-1905}
\varphi(\lambda_{n}^{\frac{2}{p-2}}\tilde{v}_n)\varphi'(\lambda_{n}^{\frac{2}{p-2}}\tilde{v}_n)\rightarrow \frac{\sqrt{2}}{2} \quad \mbox{and} \quad \varphi(\lambda_{n}^{\frac{2}{p-2}}\tilde{v}_n) \sim 2^{\frac{1}{4}}\left[\lambda_{n}^{\frac{2}{p-2}}\tilde{v}_n\right]^{\frac{1}{2}},
%\eeq
\end{equation*}
%and
%\begin{equation*}
%\varphi(\lambda_{n}^{\frac{2}{p-2}}\tilde{v}_n) \sim 2^{\frac{1}{4}}\left[\lambda_{n}^{\frac{2}{p-2}}\tilde{v}_n\right]^{\frac{1}{2}},
%\end{equation*}
%\begin{equation*}
%\beq\lab{eq:20240907-1905}
%\varphi(\lambda_{n}^{\frac{2}{p-2}}\tilde{v}_n)\varphi'(\lambda_{n}^{\frac{2}{p-2}}\tilde{v}_n)\rightarrow \frac{\sqrt{2}}{2}
%\eeq
%\end{equation*}
%and
%\begin{equation*}
%\varphi(\lambda_{n}^{\frac{2}{p-2}}\tilde{v}_n) \sim 2^{\frac{1}{4}}\left[\lambda_{n}^{\frac{2}{p-2}}\tilde{v}_n\right]^{\frac{1}{2}},
%\end{equation*}
which leads to
\begin{equation*}
%\beq\lab{eq:20240907-1906}
\lambda_{n}^{-1}\varphi(\lambda_{n}^{\frac{2}{p-2}}\tilde{v}_n(x))^{p-1}\varphi'(\lambda_{n}^{\frac{2}{p-2}}\tilde{v}_n(x))
\rightarrow 2^{\frac{p-4}{4}}\tilde{v}(x)^{\frac{p-2}{2}}.
%\eeq
\end{equation*}
Thus is, if $\Omega=\R^N$, then $\tilde{v}$ is a bounded radial positive decreasing function satisfying
\begin{equation*}
%\beq\lab{eq:20240908-1853}
-\Delta \tilde{v}=-\frac{\sqrt{2}}{2}+2^{\frac{p-4}{4}}\tilde{v}^{\frac{p-2}{2}} ~\hbox{in}~\R^N.
%\eeq
\end{equation*}
In such a case, by the following Lemma \ref{lemma:20240907-1918}, we shall see that $\tilde{v}\equiv \frac{\sqrt{2}}{2}$ in $\R^N$.

Actually, we claim that $\{\tilde{v}_n\}$ decay uniformly to zero at infinity.

We prove the claim by contradiction and assume that there exists a $\varepsilon>0$, sequences $\{\tilde{v}_n\}\subset C_{r,0}(\R^N)$ and $r_n\rightarrow +\infty$ such that $\tilde{v}_n(r_n)=\varepsilon$ and $\tilde{v}_n$ solves Eq.\eqref{eq:20240907-1840}. Set $\bar{v}_n(r):=\tilde{v}_n(r+r_n)$, then
\begin{equation*}
%\beq\lab{eq:20241108-2049}
-\left(\bar{v}''_n +\frac{(N-1)}{r+r_n}\bar{v}'_n\right)=-\varphi(\lambda_{n}^{\frac{2}{p-2}}\bar{v}_n)\varphi'(\lambda_{n}^{\frac{2}{p-2}}\bar{v}_n)
+\lambda_{n}^{-1}\varphi(\lambda_{n}^{\frac{2}{p-2}}\bar{v}_n)^{p-1}\varphi'(\lambda_{n}^{\frac{2}{p-2}}\bar{v}_n), r>-r_n.
%\eeq
\end{equation*}
Here $C_{r,0}(\R^N)$ denotes the space of continuous radial functions vanishing at $\infty$.
Passing to subsequences (still denoted by $\lambda_n$ and $\bar{v}_n$), we get that $\{\bar{v}_n\}$ converges to $\bar{v}$ in $C^{1,\alpha}(\R)\cap C_{\text{loc}}^{2}(J)$, where $J=(-\infty, R^*)$ is the support of $\bar{v}$ for some $0<R^*\leq +\infty$. Furthermore, $\bar{v}$ is a non-trivial solution of the following problem
\beq\lab{eq:20241108-2054}
\begin{cases}
-\bar{v}''=-\frac{\sqrt{2}}{2}+2^{\frac{p-4}{4}} \bar{v}^{\frac{p-2}{2}}~\hbox{in}~J,\\
\bar{v}(R^*)=\bar{v}'(R^*)=0~\hbox{if}~R^*<+\infty,
\end{cases}
\eeq
with $\bar{v}(0)=\varepsilon, \bar{v}\geq 0$ and bounded. Note that $\bar{v}_n$ is decreasing in $[-r_n,+\infty)$ and thus $\bar{v}$ is bounded and decreasing in $\R$, in particular in $J$.

{\bf Case 1: $R^*=+\infty$, i.e., $J=\R$.}
In such a case, $\bar{v}(r)$ has a limit $\bar{v}_+$ at $r=+\infty$ and a limit $\bar{v}_-$ at $r=-\infty$. In particular, $0\leq \bar{v}_+\leq \bar{v}(r)\leq \bar{v}_-<+\infty, \forall r\in \R$. Here $\bar{v}_\pm$ satisfies
$\displaystyle -\frac{\sqrt{2}}{2}+2^{\frac{p-4}{4}} \bar{v}_{\pm}^{\frac{p-2}{2}}=0.$
So $\bar{v}_+=\bar{v}_-=\frac{\sqrt{2}}{2}$ and thus $\bar{v}\equiv \frac{\sqrt{2}}{2}$, which is a contradiction to $\bar{v}(0)=\varepsilon$.

{\bf Case 2: $R^*<+\infty$.} Since also $\bar{v}$ is decreasing and bounded, it has limits $\lim_{r\rightarrow -\infty}\bar{v}(r)=\frac{\sqrt{2}}{2}$ and $\lim_{r\rightarrow -\infty}\bar{v}'(r)=0$. Multiplying by $\bar{v}'$ in both sides in the first equation of Eq.\eqref{eq:20241108-2054}, and integrating on $J$, we obtain that
\begin{equation*}
%\beq\lab{eq:20241109-1937}
\int_{-\infty}^{R^*} -\bar{v}'' \bar{v}' \mathrm{d}r=\int_{-\infty}^{R^*}\left[-\frac{\sqrt{2}}{2}+2^{\frac{p-4}{4}} \bar{v}^{\frac{p-2}{2}}\right]\bar{v}'\mathrm{d}r.
%\eeq
\end{equation*}
However,
\begin{equation*}
%\beq\lab{eq:20241109-1938}
\int_{-\infty}^{R^*} -\bar{v}'' \bar{v}' \mathrm{d}r=-\frac{1}{2} |\bar{v}'(r)|^2\Big|_{-\infty}^{R^*}=0,
%\eeq
\end{equation*}
while
\begin{equation*}
%\beq\lab{eq:20241109-1939}
\begin{aligned}
&\int_{-\infty}^{R^*}\left[-\frac{\sqrt{2}}{2}+2^{\frac{p-4}{4}} \bar{v}^{\frac{p-2}{2}}\right]\bar{v}'\mathrm{d}r
=-\frac{\sqrt{2}}{2}\bar{v}(r)\Big|_{-\infty}^{R^*} +2^{\frac{p-4}{4}} \frac{2}{p}\bar{v}(r)^{\frac{p}{2}}\Big|_{-\infty}^{R^*}\\
=&\frac{1}{2}-2^{\frac{p-4}{4}} \frac{2}{p} \left(\frac{\sqrt{2}}{2}\right)^{\frac{p}{2}}
=\frac{1}{2}-\frac{1}{p}=\frac{p-2}{2p}>0,
\end{aligned}
%\eeq
\end{equation*}
which is a contradiction again.

The claim being true, it is not possible to have $\tilde{v}\equiv \frac{\sqrt{2}}{2}$ in $\R^N$. Thus, reasoning as in the proof of Lemma \ref{lemma:20240904-0824}, we conclude that there exists some $R>0$ such that  $\Omega=B_R(0)$ and $\tilde{v}$ is a solution to the homogeneous Dirichlet-Neumann free boundary Problem \eqref{eq:20240905-1912}. Together with the uniqueness results already discussed, see the subsection \ref{uniqueness}, and Lemma \ref{lemma:20240907-1918} below, this completes the proof of Theorem \ref{th:20240905-2122}.
\hfill$\Box$

\bl\lab{lemma:20240907-1918}
If $\Omega=\R^N$, then $\tilde{v}\equiv \frac{\sqrt{2}}{2}$.
\el
\bp
If $\Omega=\R^N$, then $\tilde{v}>0$ in $\R^N$. Noting that $\tilde{v}$ is a positive radial decreasing function, $c_0:=\lim_{|x|\rightarrow +\infty}\tilde{v}(x)$ exists. In particular, one see that
$\displaystyle 0=-\frac{\sqrt{2}}{2}+2^{\frac{p-4}{4}} c_{0}^{\frac{p-2}{2}},$
which implies that $c_0=\frac{\sqrt{2}}{2}$. By the decrease of $\tilde{v}$ again, we have
$$-\Delta \tilde{v}=-\frac{\sqrt{2}}{2}+2^{\frac{p-4}{4}}\tilde{v}^{\frac{p-2}{2}}\geq 0~\hbox{in}~\R^N.$$
Set $\tilde{v}=u+\frac{\sqrt{2}}{2}$, then $u$ is a non-negative, radial, decreasing classical solution of
$$-\Delta u=h(u)\quad \mbox{in } \R^N \quad \mbox{with } h(s)=\frac{\sqrt{2}}{2}\left[(\sqrt{2}s+1)^{\frac{p-2}{2}}-1\right].$$
 Obviously, $h(s)>0$ for any $s>0$ and
$$
\lim_{s\rightarrow 0^+}\frac{h(s)}{s}=\frac{p-2}{2},\quad \lim_{s\rightarrow +\infty}\frac{h(s)}{s^{\frac{p-2}{2}}}=\frac{p-4}{4}.
$$
Taking any fixed $\gamma\in\left(1, \min\{\frac{N}{(N-2)_{+}},\frac{p-2}{2}\}\right)$, there exists some $C>0$ such that $h(s)\ge Cs^\gamma$ for all $s\ge0$. That is, $u$ is a non-negative, radial, decreasing $C^2$-function satisfying
$-\Delta u\geq C u^\gamma$ in $\R^N.$
So by \cite[Theorem 8.4]{QuittnerSouplet2007} or \cite[Theorem 2.2]{Jeanjean2024}, we conclude that $u\equiv 0$, which is equivalent to  $\tilde{v}\equiv \frac{\sqrt{2}}{2}$ in $\R^N$.
\ep
\subsection{The case $a\rightarrow +\infty$ for $1\leq N\leq 4$ as well as  $a\rightarrow a_0$ for $N\geq 5$ and proof of Theorem \ref{th:20240928-2149}}
$ $

\smallskip

{\bf Proof of Theorem \ref{th:20240928-2149}:}
Let $a^*$ be defined by Eq.\eqref{eq:20241002-2126}. Then by Corollary \ref{cro:20240904-0800},
we know that $\lambda_a\rightarrow 0$ as $a\rightarrow a^*$ and thus directly borrowing the results established in \cite[Theorem 1.3]{Genoud-Nodari-2024} one get the desired results.
\hfill$\Box$
\smallskip

Theorem \ref{th:20240928-2149} now allows to settle the cases not covered in Lemma \ref{lemma:20240903-1650}.
\bc\lab{remark:20240928-2218}
\begin{itemize}
\item[(i)]
When $N=3$ with $p \in (6,12)$ or $N=4$ with, $p \in (5,8)$, recalling the numbers $\tau_{1,p},\tau_{2,p}$ and $\tau_{3,p}$ defined in Lemma \ref{lemma:20240903-1650}, we conclude from Theorem \ref{th:20240928-2149}-(iii) that
$$ \tau_{1,p}=\|\nabla u_0\|_2^2, \quad \tau_{2,p}=V(u_0), \quad \mbox{and} \quad \tau_{3,p}=\|u_0\|_p^p.$$
Indeed, since $u_a\rightarrow u_0$ in $D_{0}^{1,2}(\R^N)\cap L^p(\R^N)$,  $\tau_{1,p}=\|\nabla u_0\|_2^2$ and $\tau_{3,p}=\|u_0\|_p^p$. On the other hand, by $P(u_a)\equiv 0$ and $P(u_0)=0$, one also see that $\tau_{2,p}=V(u_0)$.
%Thus, it holds that $u_a^2\rightarrow u_0^2$ in $D_{0}^{1,2}(\R^N)$.
\item[(ii)] For $N=3, p=6$, since the rescaling preserves the $D_{0}^{1,2}$-norm and $L^{6}$-norm, we have that
$$ \tau_{1,6}=\limsup_{a\rightarrow \infty}\|\nabla u_a\|_2^2
    = \limsup_{a\rightarrow \infty}\|\nabla \mu_{a}^{\frac{1}{2}} u_a(\mu_a \cdot)\|_2^2 =\|\nabla U\|_2^2$$
 and similarly $\tau_{3,6}=\|U\|_{6}^{6}=\tau_{1,6}.$
    Then by $P(u_a)\equiv 0$, we get that $\tau_{2,6}=\frac{1}{5}\tau_{1,6}.$
\end{itemize}
\ec
\s{Some  convergence results directly on $u_n$ as $a_n\rightarrow 0^+$ and proof of Theorem \ref{th:20250101-1347}}

As $a_n\rightarrow 0^+$, we know that
$\lambda_n\rightarrow + \infty, \, u_n(0)=\|u_n\|_\infty\rightarrow \infty, \, v_n(0)=\|v_n\|_\infty\rightarrow \infty.$
Define
$\bar{u}_n(x):= \displaystyle \lambda_{n}^{-\frac{1}{p-2}}u_n\left(\frac{x}{\lambda_{n}^{\frac{p-4}{2(p-2)}}}\right),$
then by a direct computation,
\beq\lab{eq:-7}
-\lambda_{n}^{-\frac{2}{p-2}}\Delta \bar{u}_n - \Delta(|\bar{u}_n|^2)\bar{u}_n+ \bar{u}_n
=\bar{u}_{n}^{p-1}~\hbox{in}~\R^N.
\eeq
By Lemma  \ref{lemma:20240904-1639}, and since $\|u_n\|_{\infty}=u_n(0)=\varphi(v_n(0))=\varphi(\|v_n\|_{\infty})=(2^{\frac{1}{4}}+o_n(1))\|v_n\|_{\infty}^{\frac{1}{2}},$
we know that $\{\bar{u}_n\}$ is bounded in $L^\infty(\R^N)$.
Let $\tilde{v}_n, \tilde{v}$ be defined by Theorem \ref{th:20240905-2122}.

\bl\lab{lemma:20250101-1104}
$\{V(\bar{u}_n)\}$ is bounded.
\el
\bp
Since
$$
\begin{aligned}
\bar{u}_n \nabla \bar{u}_n(x)=&\lambda_{n}^{-\frac{1}{p-2}}u_n\left(\frac{x}{\lambda_{n}^{\frac{p-4}{2(p-2)}}}\right)
\lambda_{n}^{-\frac{1}{p-2}} \, \nabla u_n \left(\frac{x}{\lambda_{n}^{\frac{p-4}{2(p-2)}}}\right) \frac{1}{\lambda_{n}^{\frac{p-4}{2(p-2)}}}\\
=&\lambda_{n}^{-\frac{2}{p-2}}\varphi\left(v_n\left(\frac{x}{\lambda_{n}^{\frac{p-4}{2(p-2)}}}\right)\right)
\varphi'\left(v_n\left(\frac{x}{\lambda_{n}^{\frac{p-4}{2(p-2)}}}\right)\right) \, \nabla v_n \left(\frac{x}{\lambda_{n}^{\frac{p-4}{2(p-2)}}}\right) \frac{1}{\lambda_{n}^{\frac{p-4}{2(p-2)}}}\\
=&\varphi\left(v_n\left(\frac{x}{\lambda_{n}^{\frac{p-4}{2(p-2)}}}\right)\right)
\varphi'\left(v_n\left(\frac{x}{\lambda_{n}^{\frac{p-4}{2(p-2)}}}\right)\right) \nabla \tilde{v}_n(x),
\end{aligned}
$$
by the fact $\varphi(s)\varphi'(s)\leq 2^{-\frac{1}{2}}, \forall s\in \mathbb{R}^+$, we have that
$\bar{u}_n |\nabla \bar{u}_n|\leq 2^{-\frac{1}{2}} |\nabla \tilde{v}_n|,$
and thus
$V(\bar{u}_n)\leq \frac{1}{2} \|\nabla \tilde{v}_n\|_2^2.$
Hence, by Theorem \ref{th:20240905-2122}, we see that $\{V(\bar{u}_n)\}$ is bounded.
\ep
By Lemma \ref{lemma:20250101-1104}, we see that $\{\bar{u}_n^2\}$ is bounded in $D_{0}^{1,2}(\R^N)$. So, we can introduce a non-negative radial function $\bar{u}$ such that, up to a subsequence, $\bar{u}_n^2\rightharpoonup \bar{u}^2$ in $D_{0}^{1,2}(\R^N)$ and $\bar{u}_n\rightarrow \bar{u}$ a.e. in $\R^N$. We also note that $\bar{u}$ is decreasing w.r.t. $r=|x|$.

\bl\lab{lemma:20250101-1152}
 Let $\Omega:=\text{supp}~\tilde{v}=B_R(0)$ be given by Theorem \ref{th:20240905-2122}. It holds $\text{supp}~\bar{u}=\Omega$.
\el
\bp
By the strict monotonicity of $\varphi$ and $\varphi(0)=0$, a direct computation shows that
\begin{align*}
\tilde{v}(x)=(1+o_n(1))\tilde{v}_n(x) = &(1+o_n(1))\lambda_{n}^{-\frac{2}{p-2}}\varphi^{-1}(\lambda_{n}^{\frac{1}{p-2}}\bar{u}_n(x))\\
%=&(1+o_n(1))\lambda_{n}^{-\frac{2}{p-2}}\varphi^{-1}(\lambda_{n}^{\frac{1}{p-2}}\bar{u}_n(x))\\
=&(1+o_n(1))\lambda_{n}^{-\frac{2}{p-2}}\varphi^{-1}\left((1+o_n(1))\bar{u}(x)\lambda_{n}^{\frac{1}{p-2}}\right).
\end{align*}
%\begin{align*}
%\tilde{v}(x)=&(1+o_n(1))\tilde{v}_n(x)\\
%=&(1+o_n(1))\lambda_{n}^{-\frac{2}{p-2}}\varphi^{-1}(\lambda_{n}^{\frac{1}{p-2}}\bar{u}_n(x))\\
%=&(1+o_n(1))\lambda_{n}^{-\frac{2}{p-2}}\varphi^{-1}\left((1+o_n(1))\bar{u}(x)\lambda_{n}^{\frac{1}{p-2}}\right).
%\end{align*}
Similarly,
$$\bar{u}(x)=(1+o_n(1))\lambda_{n}^{-\frac{1}{p-2}}\varphi\left((1+o_n(1))\tilde{v}(x)\lambda_{n}^{\frac{2}{p-2}}\right).$$
Recalling that $\lambda_n\rightarrow +\infty$ and $\varphi^{-1}(s)\sim 2^{-\frac{1}{2}}s^2, \varphi(s)\sim 2^{\frac{1}{4}}s^{\frac{1}{2}}$ near infinity, one see that $\tilde{v}(x)>0$ if and only if $\bar{u}(x)>0$.
\ep
\bl\lab{lemma:20250101-1255}
$\bar{u}_n\rightarrow \bar{u}$ in $C_{\text{loc}}^{2}(\Omega)$.
\el
\bp
Let $\Omega_1\subset\subset \Omega_2\subset\subset \Omega$. By $\inf_{x\in \Omega_2}\bar{u}>0$, we can find
some $C'_{\Omega_2}>0$ and $N_1\in \mathbb{N}$, such that,
\beq\lab{eq:20250101-1259}
0<C'_{\Omega_2}\leq \lambda_{n}^{-\frac{2}{p-2}}+2\bar{u}_n^2, \, \forall x\in \Omega_2, n\geq N_1.
\eeq
Also by $\inf_{x\in \Omega_2}\tilde{v}>0$, there exists some $C_{\Omega_2}>0$ such that
$$
v_n\left(\frac{x}{\lambda_{n}^{\frac{p-4}{2(p-2)}}}\right)
\geq C_{\Omega_2} \, \lambda_{n}^{\frac{2}{p-2}}\rightarrow +\infty~\hbox{as $n\rightarrow +\infty$ uniformly in $\Omega_2$.}
$$
Thus, using Eq.\eqref{eq:20240904-1055},
$$
%\beq\lab{eq:20250101-1304}
\varphi'\left(v_n\left(\frac{x}{\lambda_{n}^{\frac{p-4}{2(p-2)}}}\right)\right)=(1+o_n(1)) 2^{-\frac{3}{4}} \left(v_n\left(\frac{x}{\lambda_{n}^{\frac{p-4}{2(p-2)}}}\right)\right)^{-\frac{1}{2}}~\hbox{uniformly in $\Omega_2$}.
%\eeq
$$
So, by a direct calculation,
\beq\lab{eq:20250101-1305}
\begin{aligned}
\nabla \bar{u}_n(x)&= \lambda_{n}^{\frac{1}{p-2}}\varphi'\left(v_n\left(\frac{x}{\lambda_{n}^{\frac{p-4}{2(p-2)}}}\right)\right) \nabla \tilde{v}_n(x)\\
=&(1+o_n(1)) \, 2^{-\frac{3}{4}}\, \lambda_{n}^{\frac{1}{p-2}} \left(v_n\left(\frac{x}{\lambda_{n}^{\frac{p-4}{2(p-2)}}}\right)\right)^{-\frac{1}{2}}\nabla \tilde{v}_n(x)\\
=&(1+o_n(1)) \, 2^{-\frac{3}{4}} \, \frac{1}{\sqrt{\tilde{v}_n(x)}} \nabla \tilde{v}_n(x) \rightarrow \, 2^{-\frac{3}{4}} \frac{1}{\sqrt{\tilde{v}(x)}} \nabla \tilde{v}(x)~\hbox{in}~\Omega_2.
\end{aligned}
\eeq
%\beq\lab{eq:20250101-1305}
%\begin{aligned}
%\nabla \bar{u}_n(x)=&\lambda_{n}^{-\frac{1}{p-2}} \, \nabla u_n \left(\frac{x}{\lambda_{n}^{\frac{p-4}{2(p-2)}}}\right) \frac{1}{\lambda_{n}^{\frac{p-4}{2(p-2)}}}\\
%=& \lambda_{n}^{-\frac{1}{p-2}}\varphi'\left(v_n\left(\frac{x}{\lambda_{n}^{\frac{p-4}{2(p-2)}}}\right)\right) \, \nabla v_n \left(\frac{x}{\lambda_{n}^{\frac{p-4}{2(p-2)}}}\right) \frac{1}{\lambda_{n}^{\frac{p-4}{2(p-2)}}}\\
%=& \lambda_{n}^{\frac{1}{p-2}}\varphi'\left(v_n\left(\frac{x}{\lambda_{n}^{\frac{p-4}{2(p-2)}}}\right)\right) \nabla \tilde{v}_n(x)\\
%=&(1+o_n(1)) \, 2^{-\frac{3}{4}}\, \lambda_{n}^{\frac{1}{p-2}} \left(v_n\left(\frac{x}{\lambda_{n}^{\frac{p-4}{2(p-2)}}}\right)\right)^{-\frac{1}{2}}\nabla \tilde{v}_n(x)\\
%=&(1+o_n(1)) \, 2^{-\frac{3}{4}} \, \frac{1}{\sqrt{\tilde{v}_n(x)}} \nabla \tilde{v}_n(x)\\
%\rightarrow& \, 2^{-\frac{3}{4}} \frac{1}{\sqrt{\tilde{v}(x)}} \nabla \tilde{v}(x)~\hbox{in}~\Omega_2.
%\end{aligned}
%\eeq
Recalling Eq.\eqref{eq:-7}, we have that
$-(\lambda_{n}^{-\frac{2}{p-2}}+2\bar{u}_n^2) \Delta \bar{u}_n = -(1-2|\nabla \bar{u}_n|^2) \bar{u}_n+\bar{u}_{n}^{p-1}~\hbox{in}~\R^N.$
In particular,
\beq\lab{eq:20250101-1323}
-\Delta \bar{u}_n=-\frac{(1-2|\nabla \bar{u}_n|^2)}{\lambda_{n}^{-\frac{2}{p-2}}+2\bar{u}_n^2}\bar{u}_n+\frac{\bar{u}_{n}^{p-1}}{\lambda_{n}^{-\frac{2}{p-2}}+2\bar{u}_n^2}
~\hbox{in}~\Omega.
\eeq
By $\bar{u}_n\in L^\infty(\R^N)$, Eq.\eqref{eq:20250101-1259} and Eq.\eqref{eq:20250101-1305}, one check that the RHS of Eq.\eqref{eq:20250101-1323} is in $L^\infty(\Omega_2)$. So, by standard elliptic estimates, one can prove that $\bar{u}_n\rightarrow \bar{u}$ in $C^{2}(\Omega_1)$.
By the arbitrary of $\Omega_1\subset\subset \Omega_2\subset\subset \Omega$, we obtain that $\bar{u}_n\rightarrow \bar{u}$ in $C_{\text{loc}}^{2}(\Omega)$.
\ep
\bc\lab{cro:20250101-1336}
$\bar{u}(0)=2^{\frac{1}{4}}\sqrt{\tilde{v}(0)}$.
\ec
\bp
By a direct computation, using that $v_n(0) \to + \infty$,
\begin{align*}
\bar{u}(0)=&\lim_{n\rightarrow \infty} \bar{u}_n(0)
=\lim_{n\rightarrow \infty} \lambda_{n}^{-\frac{1}{p-2}}u_n(0)
=\lim_{n\rightarrow \infty} \lambda_{n}^{-\frac{1}{p-2}}\varphi(v_n(0))\\
=&\lim_{n\rightarrow \infty} \lambda_{n}^{-\frac{1}{p-2}} (1+o_n(1)) \, 2^{\frac{1}{4}} \sqrt{v_n(0)} ~\quad \hbox{by Eq.\eqref{eq3.8}}\\
=&\lim_{n\rightarrow \infty} \lambda_{n}^{-\frac{1}{p-2}} (1+o_n(1)) \, 2^{\frac{1}{4}} \lambda_{n}^{\frac{1}{p-2}} \sqrt{\tilde{v}_n(0)}
= 2^{\frac{1}{4}} \sqrt{\tilde{v}(0)}.
\end{align*}
\ep
{\bf Proof of Theorem \ref{th:20250101-1347}:}
By Lemma \ref{lemma:20250101-1255}, we can take the limit in Eq.\eqref{eq:-7} over $x\in \Omega$, and obtain that
$- \Delta(|\bar{u}|^2)\bar{u}+\bar{u}=\bar{u}^{p-1}~\hbox{in}~\Omega,$
i.e.,
$$
- \Delta(|\bar{u}|^2)+1=\bar{u}^{p-2}~\hbox{in}~\Omega.
$$
Since $\bar{u}$ is radial, by Corollary \ref{cro:20250101-1336}, we get
$\bar{u}(0)=2^{\frac{1}{4}} \sqrt{\tilde{v}(0)},~ \bar{u}'(0)=0.$
In particular, by the uniqueness of solutions of initial value problem, $\bar{u}$ is uniquely determined.
Let $w(x):=2^{-\frac{1}{2}}\bar{u}^2(x)$, a direct calculation leads to
\begin{equation*}
- \Delta w=-\frac{\sqrt{2}}{2}+2^{\frac{p-4}{4}}w^{\frac{p-2}{2}}~\hbox{in}~\Omega,~
w(0)=\tilde{v}(0), ~w'(0)=\tilde{v}'(0)=0.
\end{equation*}
So, still by the uniqueness of solutions of he IVP, we get $w=\tilde{v}$. That is, $\bar{u}=2^{\frac{1}{4}}\sqrt{\tilde{v}}$. \hfill$\Box$

\br\lab{remark:20250101-1404}
By $\bar{u}=2^{\frac{1}{4}}\sqrt{\tilde{v}}$ and $\tilde{v}=0, \frac{\partial \tilde{v}}{\partial \nu}=0$ on $\partial\Omega$, one see that  $\bar{u}=0, \bar{u}\frac{\partial \bar{u}}{\partial \nu}=0$ on $\partial\Omega$. However, we cannot assert that $\frac{\partial \bar{u}}{\partial \nu}=0$ on $\partial\Omega$ at this point.
Suppose the limit $\tilde{v}''_-(R)=\lim_{r\rightarrow R_-}\tilde{v}''(r)$ holds, then
%$$\tilde{v}''_-(R)=\lim_{r\rightarrow R_-}\tilde{v}''(r)$$ holds,
by Eq.\eqref{eq:20250101-1305} and L'Hospital's rule, formally it is expected that $\frac{\partial \bar{u}}{\partial\nu}=2^{-\frac{1}{4}} \sqrt{\tilde{v}''_-(R)}$.
The estimation of $\nabla \bar{u}_n$ on $\Omega^c$ is also unknown.
\er

\section*{Appendix: Radial symmetry of energy ground states for $a = a_0$}
\setcounter{equation}{0}
\renewcommand{\theequation}{A.\arabic{equation}}

Let $h:[0,\infty)\to [0,\infty)$ be defined by $h(s):=\varphi(s)^{p-1}\varphi'(s)$.
Define $v_0:=\varphi^{-1}(u_0)$. Then $0<v_0\in H_{\mathrm{rad}}^{1}(\mathbb{R}^N)\cap C^2(\mathbb{R}^N)$ satisfies
\begin{equation}\label{eq:20260421-1103}
-\Delta v=h(v) \quad \text{in } \mathbb{R}^N.
\end{equation}
Now take $a = a_0$, and let $u_1\in X$ be an arbitrary energy ground state.
Without loss of generality, we may assume that $u_1>0$. We note that $u_1\in X\cap C^2(\mathbb{R}^N)$ and, since $u_1^*$ is also an energy ground state, it follows that $u_1^*=u_0$.
Set $v_1:=\varphi^{-1}(u_1)$. Then $0<v_1\in H^1(\mathbb{R}^N)\cap C^2(\mathbb{R}^N)$. Define $M:=\|v_1\|_\infty$.  We see that $v_1$ also solves \eqref{eq:20260421-1103} and we observe that $ v_1 = \varphi^{-1}(u_1) \Longrightarrow v_1^* = [\varphi^{-1}(u_1)]^* = \varphi^{-1}(u_1^*)= \varphi^{-1}(u_0)=v_0$. If we manage to prove that $v_1$ is, up to translation, radial it will also be the case for $u_1$. To this aim we shall prove that $v_1= v_1^* (= v_0)$.
A standard elliptic estimate shows that
$v_1(x)\to 0$ as $|x|\to\infty$. Consequently, for every $t>0$, the function
$(v_1-t)_+\in H^1(\mathbb{R}^N)$ and has compact support.

Using the layer-cake representation, see for example \cite[Theorem 1.13]{LiLo}, and the equimeasurability of $v_1$ and $v_0$, we obtain, for any $t>0$,
\begin{equation}\label{eq:20260421-1112}
\begin{aligned}
&\int_{\mathbb{R}^N} h(v_1)(v_1-t)_+ \,\mathrm{d}x
= \int_t^\infty \frac{\mathrm{d}}{\mathrm{d}\tau}\left[h(\tau)(\tau-t)\right]\,|\{v_1>\tau\}|\,\mathrm{d}\tau \\
&= \int_t^\infty \frac{\mathrm{d}}{\mathrm{d}\tau}\left[h(\tau)(\tau-t)\right]\,|\{v_1^*>\tau\}|\,\mathrm{d}\tau
= \int_{\mathbb{R}^N} h(v_0)(v_0-t)_+ \,\mathrm{d}x.
\end{aligned}
\end{equation}
Multiplying the Eq.\eqref{eq:20260421-1103} for $v_1$ by $(v_1-t)_+$ and integrating by parts gives
\[
\int_{\mathbb{R}^N} \nabla v_1\cdot\nabla (v_1-t)_+ \,\mathrm{d}x
= \int_{\mathbb{R}^N} h(v_1)(v_1-t)_+ \,\mathrm{d}x .
\]
Performing the same computation for $v_0$ and using Eq.\eqref{eq:20260421-1112} yields
\begin{equation}\label{eq:20260421-1115}
\begin{aligned}
\int_{\{v_0>t\}} |\nabla v_0|^2 \,\mathrm{d}x
&= \int_{\mathbb{R}^N} \nabla v_0\cdot\nabla (v_0-t)_+ \,\mathrm{d}x
 = \int_{\mathbb{R}^N} h(v_0)(v_0-t)_+ \,\mathrm{d}x \\
&= \int_{\mathbb{R}^N} h(v_1)(v_1-t)_+ \,\mathrm{d}x
 = \int_{\{v_1>t\}} |\nabla v_1|^2 \,\mathrm{d}x, \qquad \forall\,0\leq t < M.
\end{aligned}
\end{equation}
Thus,
\begin{equation}\label{eq:20260421-1333}
\|\nabla (v_1-t)_+\|_2=\|\nabla (v_0-t)_+\|_2
=\|\nabla (v_1^*-t)_+\|_2=\|\nabla (v_1-t)_+^*\|_2,\quad \forall \, 0\leq t< M.
\end{equation}
Introduce the notation
$E_t:=\{v_1>t\},~E_t^*:=\{v_0>t\}$.
For any $t>0$, since $(v_1-t)_+$ has compact support and Eq.\eqref{eq:20260421-1333} holds,
%the Brothers--Ziemer result
\cite[Theorem 1.1]{BrZi} implies that there exists a vector $c(t)\in\mathbb{R}^N$ such that $E_t = E_t^*+c(t)$.
In particular, for every $\delta\in (0,M-t)$, the set $\{(v_1-t)_+>\delta\}$ shares the same center $c(t)$.

Now observe that
\begin{equation*}
%\label{eq:20260421-1441}
\big((v_1-t)_+ -(s-t)\big)_+=(v_1-s)_+, \qquad \forall \, 0<t<s<M.
\end{equation*}
Taking $\delta:=s-t\in (0, M-t)$, we obtain
$E_s=\{v_1>s\}=\{(v_1-t)_+>\delta\}$,
and therefore $c(t)=c(s)$.
Hence all the centers $c(t)$ coincide with a single vector, say $x_0\in\mathbb{R}^N$.
Consequently, every level set $\{v_1>t\}$ is a ball centered at $x_0$, which means that $v_1$ is radially symmetric about $x_0$.
%After a translation we may assume $x_0=0$.
%Because $u_1$ and $u_0$ satisfy the same equation and have the same $L^2$-norm, the uniqueness statement of Proposition~\ref{remark:20241001-1250} implies $u_1=u_0$.

%Thus, up to translation and sign change, the energy ground state for $a=a_0$ is unique and equals $u_0$.

\end{document}